\documentclass[12pt,a4paper]{amsart}
\usepackage{amssymb,amsfonts,amsthm,amsmath,graphicx,xcolor,csquotes}
\usepackage[scr]{rsfso}
\usepackage{chngcntr}
\usepackage{enumitem}

\theoremstyle{plain}
\newtheorem{theorem}{Theorem}

\numberwithin{theorem}{section}
\numberwithin{equation}{section}
\numberwithin{statement}{section}
\numberwithin{lemma}{section}
\numberwithin{definition}{section}
\numberwithin{conjecture}{section}
\numberwithin{corollary}{section}
\begin{document}
\title[Visible Point Partition Identities]{Visible Point Partition Identities for Polylogarithms, and Parametric Euler Sums}
\author{Geoffrey B Campbell}
\address{Mathematical Sciences Institute,
         The Australian National University,
         Canberra, ACT, 0200, Australia}

\email{Geoffrey.Campbell@anu.edu.au}

\thanks{Thanks to Professor Dr Henk Koppelaar, whose suggestions helped summarize, for this paper, some chapters of the author's draft book.}


\begin{abstract}
We set the scene with known values and functional relations for dilogarithms, trilogarithms and polylogarithms of various orders, along with more recent Euler sum values and multidimensional computations paying homage to the three late Professors Borwein \textit{et al.}. We then apply many of these sum values to tabulate some sixty new combinatorial identities for weighted partitions into Visible Point Vectors in 2D, 3D, 4D and 5D cases suggesting new $n$D first hyperquadrant and hyperpyramid lattice point identities.
\end{abstract}

\maketitle

\section{Early history of the Dilogarithm}\index{Dilogarithm}

The dilogarithm\index{Dilogarithm} function was first studied by Euler \cite{lE1768}\index{Euler, L.} and Landen \cite{jL1780} in the second half of the eighteenth century; and in the second half of the 20th century we had the books from Lewin \cite{lL1981,lL1982,lL1991}. So, the function known as the Euler dilogarithm\index{Dilogarithm}\index{Euler, L.}\index{Euler dilogarithm}, is defined by

\begin{equation}\label{23.01}
  Li_2(z) =  \sum_{r=1}^{\infty} \frac{z^r}{r^2}, \; \; |z|\leq 1.
\end{equation}

A version extended beyond the unit circle is given by the integral

\begin{equation}\label{23.02}
  Li_2(z) =  - \int_{0}^{z} \frac{\log(1-t)}{t} dt.
\end{equation}

For $z$ real and greater than unity the dilogarithm is complex. The function is transcendental, and is known to be irrational when $z$ is rational.

Euler easily found functional equations\index{Functional equation} and a few special values. The major functional equation\index{Functional equation} results known today are:

\begin{equation}\label{23.03}
  Li_2(z) + Li_2(-z)=  \frac{1}{2} {Li_2}^2(z),
\end{equation}

\begin{equation}\label{23.04}
  Li_2(1-z) + Li_2(1-\frac{1}{z})= - \frac{1}{2} Li_2(z^2),
\end{equation}

\begin{equation}\label{23.05}
  Li_2(z) + Li_2(1-z)=  \frac{\pi^2}{6}  - \log(z) \log(1-z),
\end{equation}

\begin{equation}\label{23.06}
  Li_2(-z) - Li_2(1-z) + \frac{1}{2} Li_2(1- z^2) =  - \frac{\pi^2}{12} - \log(z) \log(1+z).
\end{equation}

In the literature, probably the best known bivariate functional equation\index{Functional equation}\index{Functional equation!Abel's equation}\index{Functional equation!Dilogarithm} is \textbf{Abel's equation}\index{Abel's equation}, easily verified by differentiation,

\begin{equation}\label{23.07}
  Li_2\left(\frac{x}{1-x} . \frac{y}{1-y}\right) = Li_2\left(\frac{x}{1-y}\right) + Li_2\left(\frac{y}{1-x}\right) - Li_2(x) - Li_2(y) - \log(1-x) \log(1-y).
\end{equation}

This formula famously given by Abel was 20 years earlier given by Spence\index{Spence's dilogarithm equation}\index{Functional equation!Spence's equation}\index{Functional equation!Dilogarithm} in the equivalent form,

\begin{equation}\label{23.08}
  Li_2\left(\frac{x}{1-x} . \frac{y}{1-y}\right) = Li_2\left(\frac{x}{1-y}\right) + Li_2\left(\frac{y}{1-x}\right) + Li_2\left(\frac{-x}{1-x} \right)
\end{equation}
\begin{equation} \nonumber
  + Li_2\left(\frac{-y}{1-y}\right) + \frac{1}{2} {\log}^2\left(\frac{1-x}{1-y}\right).
\end{equation}

Let $L(x)$ denote the Rogers $L$-function defined in terms of the usual dilogarithm by\index{Rogers dilogarithm equation}\index{Rogers dilogarithm $L$-function}

\begin{equation}\label{23.09}
  L(x) = \frac{6}{\pi^2} \left( Li_2(x) + \frac{1}{2} \log(x) \log(1-x) \right)
       = \frac{6}{\pi^2} \left( \sum_{k=1}^{\infty} \frac{x^k}{k^2} + \frac{1}{2} \log(x) \log(1-x) \right).
\end{equation}

Then $L(x)$ satisfies the functional equation\index{Functional equation}\index{Functional equation!Rogers dilogarithm equation}

\begin{equation}\label{23.10}
  L(x) + L(y) = L(xy) + L\left(\frac{x (1-y)}{1-xy}\right) + L\left(\frac{y (1-x)}{1-xy}\right).
\end{equation}

Related to this, are the \textit{Watson identities} \cite{gW1937} defined as follows. Let $\alpha$, $-\beta$, and $-\gamma^{-1}$ be the roots of the cubic equation

 \begin{equation} \label{23.11}
 t^3+2t^2-t-1=0,
 \end{equation}

then the Rogers $L$-function satisfies\index{Functional equation!Dilogarithm cases}

\begin{equation} \label{23.12}
L(\alpha)- L(\alpha^2)	=	\frac{1}{7},	
\end{equation}
\begin{equation} \label{23.13}
L(\beta)+ \frac{1}{2} L(\beta^2)	=	\frac{5}{7},	
\end{equation}
\begin{equation} \label{23.14}
L(\gamma)+ \frac{1}{2} L(\gamma^2)	=	\frac{4}{7}.
\end{equation}

A complete list of $Li_2(x)$ which can be evaluated in closed form is given by

\begin{equation}\label{23.15}
  Li_2(-1) = - \frac{\pi^2}{12},
\end{equation}\index{Functional equation!Dilogarithm cases}
\begin{equation}\label{23.16}
  Li_2(0) = 0,
 \end{equation}\index{Functional equation!Dilogarithm cases}
\begin{equation}\label{23.17}
  Li_2(\frac{1}{2}) =  \frac{\pi^2}{12} - \frac{1}{2} (\log 2)^2,
  \end{equation}\index{Functional equation!Dilogarithm cases}
\begin{equation}\label{23.18}
  Li_2(1) = - \frac{\pi^2}{6},
  \end{equation}\index{Functional equation!Dilogarithm cases}
\begin{equation}\label{23.19}
  Li_2(2) = - \frac{\pi^2}{4} - i \pi \log2,
  \end{equation}\index{Functional equation!Dilogarithm cases}
  \begin{equation}\label{23.20}
  Li_2(-\phi) = - \frac{\pi^2}{10} - (\log \phi)^2,
  \end{equation}\index{Functional equation!Dilogarithm cases}
\begin{equation}\label{23.21}
  Li_2(-\phi^{-1}) = - \frac{\pi^2}{15} + \frac{1}{2} (\log \phi)^2,
  \end{equation}\index{Functional equation!Dilogarithm cases}
\begin{equation}\label{23.22}
  Li_2(\phi^{-2}) = \frac{\pi^2}{15} - (\log \phi)^2,
  \end{equation}\index{Functional equation!Dilogarithm cases}
\begin{equation}\label{23.23}
  Li_2(\phi^{-1}) = \frac{\pi^2}{10} - (\log \phi)^2,
\end{equation}

where $\phi$ is the golden ratio $\frac{1+\sqrt{5}}{2}$. For proof of these see Bailey et al. \cite{dB1997} and Borwein et al. \cite{jB2001}. It may seem strange that the above nine values of $Li_2(x)$ is the full known list. Also, there are many interesting identities involving the dilogarithm\index{Dilogarithm}\index{Functional equation!Dilogarithm cases} function. Ramanujan\index{Ramanujan, S.} gave the following:

\begin{equation}\label{23.24}
Li_2\left(\frac{1}{3}\right)- \frac{1}{6} Li_2\left(\frac{1}{9}\right) = \frac{\pi^2}{18} - \frac{1}{6}(\log3)^2, 	
  \end{equation}\index{Functional equation!Dilogarithm cases}
\begin{equation}\label{23.25}
Li_2\left(-\frac{1}{2}\right) + \frac{1}{6} Li_2\left(\frac{1}{9}\right)= - \frac{\pi^2}{18} + \log2 \log3 - \frac{1}{2}(\log2)^2 -\frac{1}{3}(\log3)^2  	\end{equation}\index{Functional equation!Dilogarithm cases}
\begin{equation}\label{23.26}
Li_2\left(\frac{1}{4}\right)+ \frac{1}{3} Li_2\left(\frac{1}{9}\right) = \frac{\pi^2}{18} + 2 \log2 \log3 - 2 (\log2)^2 - \frac{2}{3}(\log3)^2,  	
  \end{equation}\index{Functional equation!Dilogarithm cases}
\begin{equation}\label{23.27}
Li_2\left(-\frac{1}{3}\right)- \frac{1}{3} Li_2\left(\frac{1}{9}\right) = - \frac{\pi^2}{18} + \frac{1}{6} (\log3)^2, 	
  \end{equation}\index{Functional equation!Dilogarithm cases}
\begin{equation}\label{23.28}
Li_2\left(-\frac{1}{8}\right) + Li_2\left(\frac{1}{9}\right) = -\frac{1}{2}\left(\log\left( \frac{9}{8} \right)\right)^2.
  \end{equation}\index{Functional equation!Dilogarithm cases}

Proofs of these can be found in the Ramanujan Notebooks Vol IV\index{Ramanujan, S.} by Berndt \cite{bB1994}. In 1997, Bailey et al. \cite{dB1997} showed that

\begin{equation}\label{23.29}
 \pi^2 =36 Li_2\left(\frac{1}{2}\right) -36 Li_2\left(\frac{1}{4}\right) -12 Li_2\left(\frac{1}{8}\right) + 6 Li_2\left(\frac{1}{64}\right).
  \end{equation}\index{Functional equation!Dilogarithm cases}

Lewin \cite{lL1991} gives 67 dilogarithm\index{Dilogarithm} identities known as \textit{ladders}\index{Functional equation!Dilogarithm ladders}, and Bailey and Broadhurst \cite{dB1999} found the impressive additional \textit{17 dilogarithm} identity\index{Functional equation!Dilogarithm ladders}

 \begin{equation}\label{23.30}
 0=Li_2(\alpha_1^{-630})-2Li_2(\alpha_1^{-315})-3Li_2(\alpha_1^{-210})-10Li_2(\alpha_1^{-126})
   \end{equation}\index{Functional equation!Dilogarithm cases}
\begin{equation} \nonumber
 -7Li_2(\alpha_1^{-90})+18Li_2(\alpha_1^{-35})+84Li_2(\alpha_1^{-15})+90Li_2(\alpha_1^{-14})-4Li_2(\alpha_1^{-9})
    \end{equation}
\begin{equation} \nonumber
 +339Li_2(\alpha_1^{-8})+45Li_2(\alpha_1^{-7})+265Li_2(\alpha_1^{-6})-273Li_2(\alpha_1^{-5})-678Li_2(\alpha_1^{-4})
    \end{equation}\index{Functional equation!Dilogarithm cases}
\begin{equation} \nonumber
 -1016Li_2(\alpha_1^{-3})-744Li_2(\alpha_1^{-2})-804Li_2(\alpha_1^{-1})-22050(\log \alpha_1)^2 + 2003 \zeta(2),  	
  \end{equation}\index{Functional equation!Dilogarithm cases}

\begin{equation} \nonumber
\textmd{where} \; \alpha_1=(x^{10}+x^9-x^7-x^6-x^5-x^4-x^3+x+1)_2 \approx 1.176280818259917506544
 \end{equation}

is the largest positive root of the polynomial in \textit{Lehmer's Mahler measure problem}\index{Functional equation!Dilogarithm cases}\index{Lehmer's Mahler measure problem} and $\zeta(z)$ is the Riemann zeta function\index{Riemann zeta function}.

The closed form particular cases are all arrived at from elementary algebraic manipulation of the functional equations, and therefore will be given also as exercises. However, in these exercises clues and prompts may be useful expedients for the student.

Functional equations for the dilogarithm\index{Dilogarithm} are most easily proved by differentiating both sides. As such, the proof of the functional equations will be part of the exercises at the end of this paper.

In the literature, the Rogers $L$-function form used in (\ref{23.09}) turns out to be convenient in simplifying the trilogarithm\index{Trilogarithm} and 4th and 5th degree polylogarithm\index{Polylogarithm} functional equations.

The polylogarithm\index{Polylogarithm} occurs naturally in the Visible Point Vector identities, and the functional equations applied to these will bring relationships for weighted vector partitions.

\section{The Trilogarithm function}\index{Trilogarithm}

The trilogarithm\index{Trilogarithm} $Li_3(z)$, sometimes also denoted $L_3$, is a special case of the polylogarithm\index{Polylogarithm} $Li_n(z)$ for $n=3$. Note that the notation $Li_3(x)$ for the trilogarithm\index{Trilogarithm} is unfortunately similar to that for the logarithmic integral $Li(x)$.

Functional equations for the trilogarithm\index{Trilogarithm} include\index{Functional equation!Trilogarithm}\index{Functional equation!Trilogarithm cases}

\begin{equation}\label{23.31}
Li_3(z) + Li_3(-z) = \frac{1}{4} Li_3(z^2)
  \end{equation}\index{Functional equation!Trilogarithm}
\begin{equation}\label{23.32}
Li_3(-z) - Li_3(- z^{-1}) = -\frac{1}{6}(\log z)^3 - \frac{1}{6 \pi^2} \log z
  \end{equation}\index{Functional equation!Trilogarithm}
\begin{equation}\label{23.33}
Li_3(z) + Li_3(1-z) + Li_3 (1-z^{-1})
 = \zeta(3) + \frac{1}{6} (\log z)^3 + \frac{\pi^2}{6}\log z - \frac{1}{2}(\log z)^2 \log(1-z).
   \end{equation}\index{Functional equation!Trilogarithm}

Analytic values for $Li_3(x)$ include\index{Trilogarithm}

\begin{equation}\label{23.34}
Li_3(-1)= - \frac{3}{4} \zeta(3),
   \end{equation}\index{Trilogarithm!Values}
\begin{equation}\label{23.35}
Li_3(0)=0,
  \end{equation}\index{Trilogarithm!Values}
\begin{equation}\label{23.36}
Li_3\left(\frac{1}{2}\right)= \frac{7}{8}\zeta(3) -\frac{1}{12} \pi^2 \log 2 + \frac{1}{6} (\log 2)^3,
  \end{equation}\index{Trilogarithm!Values}
\begin{equation}\label{23.37}
Li_3(1)= \zeta(3),
  \end{equation}\index{Trilogarithm!Values}
\begin{equation}\label{23.38}
Li_3(\phi^{-2}) = \frac{4}{5} \zeta(3) - \frac{2}{15} \pi^2 \log \phi + \frac{2}{3} (\log \phi)^3,
  \end{equation}\index{Trilogarithm!Values}

where $\zeta(3)$ is Ap\'{e}ry's constant and $\phi$ is the golden ratio $\frac{1+\sqrt{5}}{2}$.

Immediate \textit{ad hoc} examples substituting the above trilogarithm values from (\ref{23.34}) to (\ref{23.38}) are available from our known derived hyperquadrant identities. The 2D, 3D, 4D and 5D hyperquadrant identities are (\ref{23.38a}), (\ref{23.38b}), (\ref{23.38c}), (\ref{23.38d}) respectively.

 If $|y|<1, |z|<1,$ and $s+t=1$, then
  \begin{equation}   \label{23.38a}
    \prod_{\substack{ (a,b)=1 \\ a,b \geq 1}} \left( \frac{1}{1-y^a z^b} \right)^{\frac{1}{a^s b^t}}
    = \exp\left\{ \left( \sum_{j=1}^{\infty} \frac{y^j}{j^s}\right) \left( \sum_{k=1}^{\infty} \frac{z^k}{k^t}\right) \right\}.
  \end{equation}

 If $|x|<1, |y|<1, |z|<1 $ and $s+t+u=1$, then
  \begin{equation}   \label{23.38b}
    \prod_{\substack{ (a,b,c)=1 \\ a,b,c \geq 1}} \left( \frac{1}{1-x^a y^b z^c} \right)^{\frac{1}{a^s b^t c^u}}
    = \exp\left\{ \left( \sum_{i=1}^{\infty} \frac{x^i}{i^s}\right) \left( \sum_{j=1}^{\infty} \frac{y^j}{j^t}\right)
    \left( \sum_{k=1}^{\infty} \frac{z^k}{k^u}\right)\right\}.
  \end{equation}

 If $|w|<1, |x|<1, |y|<1, |z|<1 $ and $r+s+t+u=1$, then
  \begin{equation}   \label{23.38c}
    \prod_{\substack{ (a,b,c,d)=1 \\ a,b,c,d \geq 1}} \left( \frac{1}{1-w^a x^b y^c z^d} \right)^{\frac{1}{a^r b^s c^t d^u}}
    = \exp\left\{ \left( \sum_{h=1}^{\infty} \frac{w^h}{h^r}\right) \left( \sum_{i=1}^{\infty} \frac{x^i}{i^s}\right)
    \left( \sum_{j=1}^{\infty} \frac{y^j}{j^t}\right) \left( \sum_{k=1}^{\infty} \frac{z^k}{k^u}\right)\right\}.
  \end{equation}

 If $|v|<1, |w|<1, |x|<1, |y|<1, |z|<1 $ and $q+r+s+t+u=1$, then
  \begin{equation}   \label{23.38d}
    \prod_{\substack{ (a,b,c,d,e)=1 \\ a,b,c,d,e \geq 1}} \left( \frac{1}{1-v^a w^b x^c y^d z^e} \right)^{\frac{1}{a^q b^r c^s d^t e^u}}
    \end{equation}
  \begin{equation} \nonumber
  = \exp\left\{ \left( \sum_{g=1}^{\infty} \frac{v^g}{g^q}\right) \left( \sum_{h=1}^{\infty} \frac{w^h}{h^r}\right) \left( \sum_{i=1}^{\infty} \frac{x^i}{i^s}\right)
    \left( \sum_{j=1}^{\infty} \frac{y^j}{j^t}\right) \left( \sum_{k=1}^{\infty} \frac{z^k}{k^u}\right)\right\}.
  \end{equation}

\subsection{2D Hyperquadrant VPV Cases}
Equation (\ref{23.38a}) has the following three example cases using our known sums for the trilogarithm and $Li_2(z)$.
 If  $|z|<1$ then
  \begin{equation}   \label{23.38a1}
    \prod_{\substack{ (a,b)=1 \\ a,b \geq 1}} \left( 1-(-1)^a z^b \right)^{\frac{b^{2}}{a^3}}
    = \exp\left\{ \frac{3z(1+z)\zeta(3)}{4(1-z)^3} \right\},
  \end{equation}
 by (\ref{23.34}) and $Li_{-2}(z) = \sum_{k=1}^{\infty} k^2 z^k = (z (1+z))/(1-z)^3$ when $|z|<1$.

 If  $|z|<1$ then
  \begin{equation}   \label{23.38a2}
    \prod_{\substack{ (a,b)=1 \\ a,b \geq 1}} \left( 1- \frac{z^b}{2^a} \right)^{\frac{b^{2}}{a^3}}
    = \exp\left\{\left(\frac{7}{8}\zeta(3) -\frac{1}{12} \pi^2 \log 2 + \frac{1}{6} (\log 2)^3 \right)\left(\frac{z (1+z)}{(1-z)^3}\right) \right\},
  \end{equation}
 by (\ref{23.36}) and our value for $Li_{-2}(z)$.

 If  $|z|<1$ then
  \begin{equation}   \label{23.38a3}
    \prod_{\substack{ (a,b)=1 \\ a,b \geq 1}} \left( 1- \frac{z^b}{\phi^{2a}} \right)^{\frac{b^{2}}{a^3}}
    = \exp\left\{\left(\frac{4}{5} \zeta(3) - \frac{2}{15} \pi^2 \log \phi + \frac{2}{3} (\log \phi)^3 \right)\left(\frac{z (1+z)}{(1-z)^3}\right) \right\},
  \end{equation}
 with $\phi$ the golden ratio; substituting (\ref{23.38}) and the value of $Li_{-2}(z)$.

A good way to summarize what can be deduced from the 2D VPV first quadrant identity from known dilogarithm and trilogarithm sums is to use the two identities:
 If $|y|<1, |z|<1,$ then
  \begin{equation}   \label{23.38a4}
    \prod_{\substack{ (a,b)=1 \\ a,b \geq 1}} \left( \frac{1}{1-y^a z^b} \right)^{\frac{b}{a^2}}
    = \exp\left\{ \frac{z Li_2(y)}{(1-z)^2} \right\}
  \end{equation}
 and
  \begin{equation}   \label{23.38a5}
    \prod_{\substack{ (a,b)=1 \\ a,b \geq 1}} \left( \frac{1}{1-y^a z^b} \right)^{\frac{b^2}{a^3}}
    = \exp\left\{ \frac{z(1+z) Li_3(y)}{(1-z)^3} \right\}.
  \end{equation}

\subsection{3D Hyperquadrant VPV Cases}
Next, the 3D first hyperquadrant VPV identity is:
 If $|x|<1, |y|<1, |z|<1 $ and $s+t+u=1$, then
  \begin{equation}   \label{23.38b1}
    \prod_{\substack{ (a,b,c)=1 \\ a,b,c \geq 1}} \left( \frac{1}{1-x^a y^b z^c} \right)^{\frac{1}{a^s b^t c^u}}
    = \exp\left\{ \left( \sum_{i=1}^{\infty} \frac{x^i}{i^s}\right) \left( \sum_{j=1}^{\infty} \frac{y^j}{j^t}\right)
    \left( \sum_{k=1}^{\infty} \frac{z^k}{k^u}\right)\right\}.
  \end{equation}
Using the nine known values of $L_2(z)$, and the five known values of $L_3(z)$, applied to (\ref{23.38b1}) we get cases of the kinds given from (\ref{23.38b2}) to (\ref{23.38b7}) here:

If $|x|<1, |y|<1, |z|<1 $ then
  \begin{equation}   \label{23.38b2}
    \prod_{\substack{ (a,b,c)=1 \\ a,b,c \geq 1}} \left( \frac{1}{1-x^a y^b z^c} \right)^{\frac{c}{a b}}
    = \exp\left\{Li_1(x) Li_1(y) Li_{-1}(z)\right\} = \left( 1-x  \right)^{\frac{z \log(1-y)}{(1-z)^2}},
  \end{equation}

  \bigskip

  \begin{equation}   \label{23.38b3}
    \prod_{\substack{ (a,b,c)=1 \\ a,b,c \geq 1}} \left( \frac{1}{1-x^a y^b z^c} \right)^{\frac{c^2}{a b^2}}
    = \exp\left\{Li_1(x) Li_2(y) Li_{-2}(z)\right\}
   \end{equation}
  \begin{equation}   \nonumber
    = \left( 1-x  \right)^{\frac{z(1+z) Li_2(y)}{(1-z)^3}},
  \end{equation}

  \bigskip

  \begin{equation}   \label{23.38b4}
    \prod_{\substack{ (a,b,c)=1 \\ a,b,c \geq 1}} \left( \frac{1}{1-x^a y^b z^c} \right)^{\frac{c^3}{a b^3}}
    = \exp\left\{Li_1(x) Li_3(y) Li_{-3}(z)\right\}
  \end{equation}
  \begin{equation}   \nonumber
    = \left( \frac{1}{1-x} \right)^{\frac{z(1+4z+z^2) Li_3(y)}{(1-z)^4}},
  \end{equation}

  \bigskip

  \begin{equation}   \label{23.38b5}
    \prod_{\substack{ (a,b,c)=1 \\ a,b,c \geq 1}} \left( \frac{1}{1-x^a y^b z^c} \right)^{\frac{c^3}{a^2 b^2}}
    = \exp\left\{Li_2(x) Li_2(y) Li_{-3}(z)\right\}
  \end{equation}
  \begin{equation}   \nonumber
    = exp\left\{ \frac{z(1+4z+z^2) Li_2(x) Li_2(y)}{(1-z)^4} \right\},
  \end{equation}

  \bigskip

  \begin{equation}   \label{23.38b6}
    \prod_{\substack{ (a,b,c)=1 \\ a,b,c \geq 1}} \left( \frac{1}{1-x^a y^b z^c} \right)^{\frac{c^4}{a^2 b^3}}
    = \exp\left\{ Li_2(x) Li_3(y) Li_{-4}(z) \right\}
  \end{equation}
  \begin{equation}   \nonumber
    = exp\left\{ \frac{z(1+11z+11z^2+z^3) Li_2(x) Li_3(y)}{(1-z)^5} \right\},
  \end{equation}

  \bigskip

  \begin{equation}   \label{23.38b7}
    \prod_{\substack{ (a,b,c)=1 \\ a,b,c \geq 1}} \left( \frac{1}{1-x^a y^b z^c} \right)^{\frac{c^5}{a^3 b^3}}
    = \exp\left\{Li_3(x) Li_3(y) Li_{-5}(z)\right\}
  \end{equation}
  \begin{equation}   \nonumber
    = exp\left\{ \frac{z(1+26z+66z^2+26z^3+z^4) Li_3(x) Li_3(y)}{(1-z)^6} \right\}.
  \end{equation}

\subsection{4D Hyperquadrant VPV Cases}
 If $|w|<1, |x|<1, |y|<1, |z|<1 $ and $r+s+t+u=1$, then
  \begin{equation}   \label{23.38c1}
    \prod_{\substack{ (a,b,c,d)=1 \\ a,b,c,d \geq 1}} \left( \frac{1}{1-w^a x^b y^c z^d} \right)^{\frac{1}{a^r b^s c^t d^u}}
    = \exp\left\{ Li_r(w) Li_{s}(x) Li_{t}(y) Li_u(z)\right\}.
  \end{equation}

We can apply all of our known dilogarithm\index{Dilogarithm} and trilogarithm\index{Trilogarithm} sums as well as use the facts $Li_1(z)= -\log(1-z)$, $Li_0(z)= z/(1-z)$, with also the known values
\begin{equation} \nonumber
  Li_{-n}(z) = \sum_{k=0}^{n}k! S(n+1,k+1) \left( \frac{z}{1-z} \right)^{k+1} \quad (n=1,2,3,\ldots),
\end{equation}
where $S(n,k)$ are Stirling numbers of the second kind.\index{Stirling numbers of the second kind}
We see that for (\ref{23.38c1}), most of the possible cases with closed form RHS's of the resulting identity case are given by the table:

\bigskip
   \begin{table}[ht]
  \centering
  \caption{Cases of $\mathbf{r+s+t+u=1}$ related to the 4D Hyperquadrant VPV} \label{Table 28}
\begin{equation} \nonumber
\begin{array}{|c|c|c|c|c|} \hline
  r & s & t & u & \log\textmd{(RHS of 4D identity)}    \\ \hline
  3 & -2 & -3 & 3 & Li_3(w) Li_{-2}(x) Li_{-3}(y) Li_3(z)  \\ \hline
  2 & -2 & -2 & 3 & Li_2(w) Li_{-2}(x) Li_{-2}(y) Li_3(z)  \\ \hline
  1 & -2 & -1 & 3 & Li_1(w) Li_{-2}(x) Li_{-1}(y) Li_3(z)  \\ \hline
  3 & -1 & -3 & 2 & Li_3(w) Li_{-1}(x) Li_{-3}(y) Li_2(z)  \\ \hline
  2 & -1 & -2 & 2 & Li_2(w) Li_{-1}(x) Li_{-2}(y) Li_2(z)  \\ \hline
  1 & -1 & -1 & 2 & Li_1(w) Li_{-1}(x) Li_{-1}(y) Li_2(z)  \\ \hline
  0 & -1 & -1 & 3 & Li_0(w) Li_{-1}(x) Li_{-1}(y) Li_{3}(z)  \\ \hline
  0 & 0 & -1 & 2 & Li_0(w) Li_0(x) Li_{-1}(y) Li_{2}(z)  \\ \hline
  0 & 0 & 0 & 1 & Li_0(w) Li_0(x) Li_0(y) Li_1(z)  \\ \hline
  0 & 1 & 1 & -1 & Li_0(w) Li_1(x) Li_1(y) Li_{-1}(z)  \\ \hline
  1 & 1 & 1 & -2 & Li_1(w) Li_1(x) Li_1(y) Li_{-2}(z)  \\ \hline
  1 & 1 & 2 & -3 & Li_1(w) Li_1(x) Li_2(y) Li_{-3}(z)  \\ \hline
  1 & 2 & 2 & -4 & Li_1(w) Li_2(x) Li_2(y) Li_{-4}(z)  \\ \hline
  2 & 2 & 2 & -5 & Li_2(w) Li_2(x) Li_2(y) Li_{-5}(z)  \\ \hline
  1 & 1 & 3 & -4 & Li_1(w) Li_1(x) Li_3(y) Li_{-4}(z)  \\ \hline
  1 & 2 & 3 & -5 & Li_1(w) Li_2(x) Li_3(y) Li_{-5}(z)  \\ \hline
  2 & 2 & 3 & -6 & Li_2(w) Li_2(x) Li_3(y) Li_{-6}(z)  \\ \hline
  2 & 3 & 3 & -7 & Li_2(w) Li_3(x) Li_3(y) Li_{-7}(z)  \\ \hline
  3 & 3 & 3 & -8 & Li_3(w) Li_3(x) Li_3(y) Li_{-8}(z)  \\ \hline
\end{array}
\end{equation}
\end{table}
The above table has values of $r$, $s$, $t$, and $u$ each less than or equal to 3, where $r+s+t+u=1$; and the resulting cases of  (\ref{23.38c1}).

\subsection{5D Hyperquadrant VPV Cases}

  \begin{table}[!t]
  \centering
  \caption{Cases of $\mathbf{q+r+s+t+u=1}$ related to the 5D Hyperquadrant VPV} \label{Table 29}
\begin{equation} \nonumber
\begin{array}{|c|c|c|c|c|c|} \hline
  q & r & s & t & u & \log\textmd{(RHS of 5D identity (\ref{23.38d1}))}    \\ \hline
  0 & 3 & -2 & -3 & 3 & Li_0(v) Li_3(w) Li_{-2}(x) Li_{-3}(y) Li_3(z)  \\ \hline
  0 & 2 & -2 & -2 & 3 & Li_0(v) Li_2(w) Li_{-2}(x) Li_{-2}(y) Li_3(z)  \\ \hline
  0 & 1 & -2 & -1 & 3 & Li_0(v) Li_1(w) Li_{-2}(x) Li_{-1}(y) Li_3(z)  \\ \hline
  0 & 0 & -1 & -1 & 3 & Li_0(v) Li_0(w) Li_{-1}(x) Li_{-1}(y) Li_{3}(z)  \\ \hline
  0 & 3 & -1 & -3 & 2 & Li_0(v) Li_3(w) Li_{-1}(x) Li_{-3}(y) Li_2(z)  \\ \hline
  0 & 2 & -1 & -2 & 2 & Li_0(v) Li_2(w) Li_{-1}(x) Li_{-2}(y) Li_2(z)  \\ \hline
  0 & 1 & -1 & -1 & 2 & Li_0(v) Li_1(w) Li_{-1}(x) Li_{-1}(y) Li_2(z)  \\ \hline
  0 & 0 & 0 & -1 & 2 & Li_0(v) Li_0(w) Li_0(x) Li_{-1}(y) Li_{2}(z)  \\ \hline
  0 & 0 & 0 & 0 & 1 & Li_0(v) Li_0(w) Li_0(x) Li_0(y) Li_1(z)  \\ \hline
  0 & 0 & 1 & 1 & -1 & Li_0(v) Li_0(w) Li_1(x) Li_1(y) Li_{-1}(z)  \\ \hline
  0 & 1 & 1 & 1 & -2 & Li_0(v) Li_1(w) Li_1(x) Li_1(y) Li_{-2}(z)  \\ \hline
  0 & 1 & 1 & 2 & -3 & Li_0(v) Li_1(w) Li_1(x) Li_2(y) Li_{-3}(z)  \\ \hline
  0 & 1 & 2 & 2 & -4 & Li_0(v) Li_1(w) Li_2(x) Li_2(y) Li_{-4}(z)  \\ \hline
  0 & 2 & 2 & 2 & -5 & Li_0(v) Li_2(w) Li_2(x) Li_2(y) Li_{-5}(z)  \\ \hline
  0 & 1 & 1 & 3 & -4 & Li_0(v) Li_1(w) Li_1(x) Li_3(y) Li_{-4}(z)  \\ \hline
  0 & 1 & 2 & 3 & -5 & Li_0(v) Li_1(w) Li_2(x) Li_3(y) Li_{-5}(z)  \\ \hline
  0 & 2 & 2 & 3 & -6 & Li_0(v) Li_2(w) Li_2(x) Li_3(y) Li_{-6}(z)  \\ \hline
  0 & 2 & 3 & 3 & -7 & Li_0(v) Li_2(w) Li_3(x) Li_3(y) Li_{-7}(z)  \\ \hline
  0 & 3 & 3 & 3 & -8 & Li_0(v) Li_3(w) Li_3(x) Li_3(y) Li_{-8}(z)  \\ \hline
  1 & 0 & 1 & 1 & -2 & Li_0(v) Li_0(w) Li_1(x) Li_1(y) Li_{-1}(z)  \\ \hline
  1 & 1 & 1 & 1 & -3 & Li_0(v) Li_1(w) Li_1(x) Li_1(y) Li_{-2}(z)  \\ \hline
  1 & 1 & 1 & 2 & -4 & Li_0(v) Li_1(w) Li_1(x) Li_2(y) Li_{-3}(z)  \\ \hline
  1 & 1 & 2 & 2 & -5 & Li_0(v) Li_1(w) Li_2(x) Li_2(y) Li_{-4}(z)  \\ \hline
  1 & 2 & 2 & 2 & -6 & Li_0(v) Li_2(w) Li_2(x) Li_2(y) Li_{-5}(z)  \\ \hline
  1 & 1 & 1 & 3 & -5 & Li_0(v) Li_1(w) Li_1(x) Li_3(y) Li_{-4}(z)  \\ \hline
  1 & 1 & 2 & 3 & -6 & Li_0(v) Li_1(w) Li_2(x) Li_3(y) Li_{-5}(z)  \\ \hline
  1 & 2 & 2 & 3 & -7 & Li_0(v) Li_2(w) Li_2(x) Li_3(y) Li_{-6}(z)  \\ \hline
  1 & 2 & 3 & 3 & -8 & Li_0(v) Li_2(w) Li_3(x) Li_3(y) Li_{-7}(z)  \\ \hline
  1 & 3 & 3 & 3 & -9 & Li_0(v) Li_3(w) Li_3(x) Li_3(y) Li_{-8}(z)  \\ \hline
\end{array}
\end{equation}
\end{table}

Next we consider the range of closed form RHS versions of the 5D hyperquadrant identity (\ref{23.38d1}). Table \ref{Table 29} of finitely evaluable cases would include all integers $q$,$r$,$s$,$t$,$u$ less than or equal to 3 such that $\mathbf{q+r+s+t+u=1}$.
The identity related to this is as follows.
 If $|v|<1, |w|<1, |x|<1, |y|<1, |z|<1 $ and $q+r+s+t+u=1$, then
  \begin{equation}   \label{23.38d1}
    \prod_{\substack{ (a,b,c,d,e)=1 \\ a,b,c,d,e \geq 1}} \left( \frac{1}{1-v^a w^b x^c y^d z^e} \right)^{\frac{1}{a^q b^r c^s d^t e^u}}
    \end{equation}
  \begin{equation} \nonumber
  = \exp\left\{ Li_q(v) Li_r(w) Li_{s}(x) Li_{t}(y) Li_u(z)\right\}.
  \end{equation}

  Refer to table \ref{Table 29}.

\newpage
\subsection{More Trilogarithm equations}
In 1997 Bailey et al., \cite{dB1997} showed that\index{Trilogarithm}

\begin{equation} \label{23.39}
\frac{35}{2} \zeta(3) - \pi^2 \log 2
= 36 Li_3\left(\frac{1}{2}\right) - 18 Li_3\left(\frac{1}{4}\right)- 4 Li_3\left(\frac{1}{8}\right) + Li_3\left(\frac{1}{64}\right),
   \end{equation}\index{Trilogarithm!Identities}
\begin{equation} \label{23.40}
 2 (\log 2)^3 - 7 \zeta(3)
 =  -24 Li_3\left(\frac{1}{2}\right) + 18 Li_3\left(\frac{1}{4}\right) + 4 Li_3\left(\frac{1}{8}\right) - Li_3\left(\frac{1}{64}\right),
   \end{equation}\index{Trilogarithm!Identities}
\begin{equation} \label{23.41}
10 (\log 2)^3 -2 \pi^2 \log 2
= -48 Li_3\left(\frac{1}{2}\right) + 54 Li_3\left(\frac{1}{4}\right) +12 Li_3\left(\frac{1}{8}\right) - 3 Li_3\left(\frac{1}{64}\right).
\end{equation}\index{Trilogarithm!Identities}

These relationships are applicable to transforms of VPV identities such as those given in this paper.

\section{The Polylogarithm function}\index{Polylogarithm}

The polylogarithm\index{Polylogarithm} is the obvious extension of the dilogarithm\index{Dilogarithm} and trilogarithm\index{Trilogarithm}. The $n$th order polylogarithm\index{Polylogarithm} is defined by

\begin{equation}\label{23.42}
  Li_n(z) =  \sum_{r=1}^{\infty} \frac{z^r}{r^n}, \; \; |z|\leq 1.
\end{equation}

A version extended beyond the unit circle is given by the integral

\begin{equation}\label{23.43}
  Li_n(z) =  \int_{0}^{z} \frac{Li_{n-1}(t)}{t} dt.
\end{equation}

For $z$ real and greater than unity the polylogarithm\index{Polylogarithm} is complex.

Research on this function extends back over 200 years, but in the past 60 years, the main literature for it was written by L Lewin \cite{lL1991}. However, with the rise of experimental computer science calculations by notably J Borwein \cite{dB2006} to \cite{jB2003} and his colleagues, both Euler sums\index{Euler, L.} and MTW sums and generalized polylogarithm\index{Polylogarithm} cases and functional equations have been found. The work of Borwein has been important for the fact that his new approaches for cases of polylogarithms\index{Polylogarithm} can find direct application to multdimensional Visible Point Vector identities of the kinds considered in the present volume. Since Jon Borwein was truly a pioneer of the calculation and discovery methods it is worth showing his image. Incidentally, since Borwein's passing in 2016, there has been a wealth of summations of Euler series\index{Euler, L.} and polylogarithm cases published by Anthony Sofo \index{Sofo, A.} \cite{aS2008} to \cite{aS2019}, many of which are applicable to Visible Point Vector Identities.

So prolific is Prof Borwein's output, \textit{Microsoft Academic} has listed him as number 21 on their list of top mathematics authors of all time. He has also produced over a dozen books, including a monograph on convex functions which was a \textit{Choice 2011} Outstanding Academic Book, and over 450 refereed articles and proceedings.

 \begin{figure} [ht!]
\centering
    \includegraphics[width=8cm,angle=0,height=8cm]{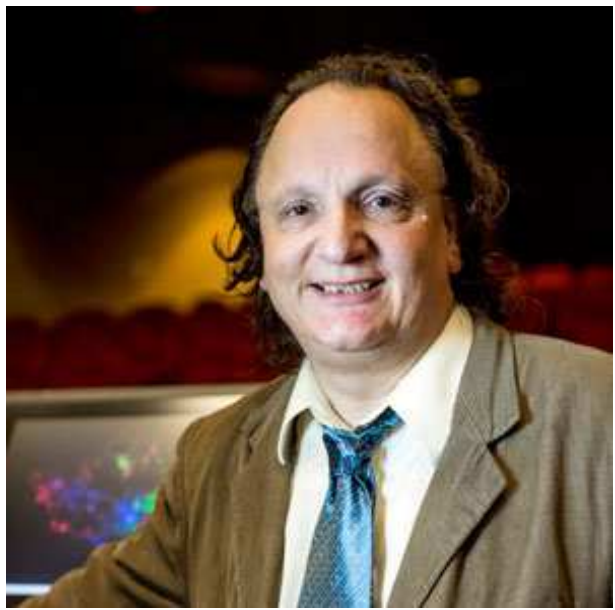}
  \caption[Picture of Professor Jonathan M. Borwein circa 2015]{Professor Jonathan M. Borwein.} \label{Fig20}
\end{figure}


\section{Mordell-Tornheim-Witten ensembles}\index{Mordell-Tornheim-Witten (MTW) zeta function}

The Borweins contributed to experimentally calculated values of the multidimensional Mordell-Tornheim-Witten (MTW) zeta function.\index{Mordell-Tornheim-Witten (MTW) zeta function}

\begin{equation} \label{23.44}
  \omega(s_1, ... , s_{K+1}) = \sum_{m_1,...,m_K >0} \frac{1}{{m_1}^{s_1} {m_2}^{s_2} \ldots {m_K}^{s_{K}}(m_1 + m_2 + \ldots +m_{K})^{s_{K+1}}}
\end{equation}

This MTW function has known relations going back 70 years, but remains mysterious with respect to many combinatorial
phenomena, especially when we contemplate derivatives with respect to the $s_i$ parameters.
We refer to $K + 1$ as the \textit{depth} and $\sum_{j=1}^{k+1} s_j$ as the \textit{weight} of $\omega$.

We can get a more intuitive grasp of MTWs by writing the simplest ones to see what is going on with them. Hence we have

\begin{equation} \label{23.44.1}
  \omega(s_1, s_{2}) = \sum_{m_1 >0} \frac{1}{{m_1}^{s_1}(m_1)^{s_{2}}},
\end{equation}\index{Mordell-Tornheim-Witten (MTW) zeta function}
\begin{equation} \label{23.44.2}
  \omega(s_1, s_2 , s_{3}) = \sum_{m_1,\, m_2 >0} \frac{1}{{m_1}^{s_1} {m_2}^{s_2}(m_1 + m_2)^{s_{3}}},
\end{equation}\index{Mordell-Tornheim-Witten (MTW) zeta function}
\begin{equation} \label{23.44.3}
  \omega(s_1, s_2 , s_{3}, s_{4}) = \sum_{m_1,\, m_2,\, m_3 >0} \frac{1}{{m_1}^{s_1} {m_2}^{s_2}{m_3}^{s_3}(m_1 + m_2+m_3)^{s_{4}}},
\end{equation}\index{Mordell-Tornheim-Witten (MTW) zeta function}

There is a hyperpyramid lattice VPV identity first given in Campbell\index{Campbell, G.B.} \cite{gC1994a}, which is evidently including MTW zeta functions.\index{Mordell-Tornheim-Witten (MTW) zeta function} In fact we have

 \begin{equation} \label{23.44a}
\prod_{\substack{(a_1, a_2, \dots, a_n)=1\\ a_1, \dots,
a_{n-1}<a_n \\
a_1, \dots, a_{n-1} \geq 0\\
a_n \geq 1}} \left(\frac{1}{1-x^{a_{1}}_1 x_2^{a_{2}} \dots
x_n^{a_{n}}}\right)^{\frac{1}{a_1^{b_{1}} a_2^{b_{2}} \dots
a_n^{b_{n}}}} =\exp \left\{ \sum^\infty_{k=1} \prod^{n-1}_{i=1}
\left(\sum^{k-1}_{j=1} \frac{x^j_i}{j^{b_{i}}}\right)
\frac{x^k_n}{k^{b_{n}}}\right\}
 \end{equation}

where for  $i = 1, 2, ..., n$;  $|x_i| < 1$,  $b_i$  is a complex
number with $\sum^n_{i=1} b_i=1$.

The principle underpinning the VPV identities is that we can apply lattice summations where the shape of the lattice in nD space is radial from the origin, thus defining continuous $n$-variate functions that sum exactly to elementary functions. This concept has not fully been utilized in the present literature.

After two decades of \textit{Generalized Euler Sum}\index{Euler, L.} and \textit{Polylogarithm Sum}\index{Euler sum} research by Bailey, Bradley and the three Borweins setting the scenes (see Bailey-Borwein(David, Peter and Jon)-Bradley-Broadhurst-Plouffe \cite{dB1995}, \cite{dB1997}, \cite{dB1999}, \cite{jB1998}, \cite{jB2001} together in varied combinations), the 2016 work Bailey et. al. \cite{dB2016} introduced and discussed a novel generalized MTW ensemble
zeta function for positive integers $M, \, N$ and nonnegative integers $s_i, \, t_j$ with constraints
$M \geq N \geq 1$ together with a polylogarithm-integral representation\index{Polylogarithm}:

\begin{equation}\label{23.45}
  \omega(s_1, ... , s_{M} \mid t_1, ... , t_{N} ) = \sum_{\substack{ m_1, ... , m_{M}, \,n_1, ... , n_{N}>0 \\
  \sum_{j=1}^{M} m_j = \sum_{k=1}^{N} n_k}} \prod_{j=1}^{M} \frac{1}{{m_j}^{s_j}} \prod_{k=1}^{N} \; \frac{1}{{n_k}^{t_k}}
\end{equation}

\begin{equation}\label{23.46}
  = \frac{1}{2\pi} \int_{0}^{2\pi} \prod_{j=1}^{M} Li_{s_j} \left( e^{i\theta} \right) \prod_{k=1}^{N} Li_{t_k} \left( e^{-i\theta} \right) d\theta.
\end{equation}

Here the \textit{polylogarithm of order s}\index{Polylogarithm} denotes $Li_s(z) := \sum_{n\geq 0} \frac{z^n}{n^s}$ and its analytic extensions
[26] and the (complex) number $s$ is its order.
Note that if parameters are zero, there are convergence issues with this integral representation. One may use principal-value calculus, or an alternative representation such as (\ref{23.48}) below. When $N = 1$ the representation ((\ref{23.47}) devolves to the classic MTW form, in that

\begin{equation} \label{23.47}
  \omega(s_1, ... , s_{M+1}) = \omega(s_1, ... , s_M \, \mid \, s_{M+1}).
\end{equation}

Because they represent examples of lattice points distributed radially from the origin in multidimensional space, these Euler Sum/Polylogarithm Sum/MTW Sum combinations\index{Euler, L.} are all applicable to the VPV identities, and can be adapted to the theorems on $n$D hyperquadrant, hyperpyramid and skewed hyperpyramid summations. Collecting these Summations with the cases arising from the more miscellaneous Euler Sum\index{Euler, L.}\index{Euler sum} evaluations given in the 4o or more applicable papers by Sofo \index{Sofo, A.} from 2008 to circa 2020, we see a vista opening wherein the numerous cases and theories can be applied to the VPV identities, and it seems hundreds of examples could be easily written from them. The way forward through all this may be to seek applications and physical models that represent particular examples where VPV identities arise.

So, returning to our summations, we need a wider \textit{MTW ensemble} involving outer derivatives, according to

\begin{equation} \label{23.48}
  \omega \left(\substack{ s_1, ... , s_{M}, \, | \,t_1, ... , t_{N} \\
                          d_1, ... , d_{M}, \, | \,e_1, ... , e_{N}} \right)
  = \sum_{\substack{ m_1, ... , m_{M}, \,n_1, ... , n_{N}>0 \\
  \sum_{j=1}^{M} m_j = \sum_{k=1}^{N} n_k}} \prod_{j=1}^{M} \frac{(-\log m_j)^{d^j}}{{m_j}^{s_j}} \prod_{k=1}^{N} \; \frac{(-\log n_k)^{e^k}}{{n_k}^{t_k}}
\end{equation}

\begin{equation}\label{23.49}
  = \frac{1}{2\pi} \int_{0}^{2\pi} \prod_{j=1}^{M} {Li_{s_j}}^{(d_j)} \left( e^{i\theta} \right) \prod_{k=1}^{N} {Li_{t_k}}^{(e_k)} \left( e^{-i\theta} \right) d\theta,
\end{equation}

where the $s$-th outer derivative of a polylogarithm\index{Polylogarithm} is denoted ${Li_s}^{(d)}(z) = \left( \frac{\partial}{\partial z} \right)^{(d)} Li_s(z)$.

\section{Finite Euler Sums}\index{Finite Euler sums}\index{Euler, L.}

Consider the following well known summations:

\begin{equation} \label{23.50}
\sum_{k=1}^{n} k^1 = \frac{n}{2} + \frac{n^2}{2},
\end{equation}\index{Finite Euler sums}
\begin{equation} \label{23.51}
\sum_{k=1}^{n} k^2 = \frac{n}{6} + \frac{n^2}{2} + \frac{n^3}{3},
\end{equation}\index{Finite Euler sums}
\begin{equation} \label{23.52}
\sum_{k=1}^{n} k^3 = \frac{n^2}{4} + \frac{n^3}{3} + \frac{n^4}{4},
\end{equation}\index{Finite Euler sums}
\begin{equation} \label{23.53}
\sum_{k=1}^{n} k^4 = \frac{-n}{30} + \frac{n^3}{3} + \frac{n^4}{2} + \frac{n^5}{5},
\end{equation}\index{Finite Euler sums}

and their generalizations:

\begin{equation} \label{23.54}
\sum_{k=1}^{n} k^1 z^k = \frac{z - (1 + n) z^{n+1} + n z^{n+2}}{(1-z)^2},
\end{equation}\index{Finite Euler sums}
\begin{equation} \label{23.55}
\sum_{k=1}^{n} k^2 z^k
= \frac{ - n^2 z^{n+3} + (2n^2+2n-1)z^{n+2} -(n^2+2n+1) z^{n+1} + z^2 + z}{(1-z)^3},
\end{equation}\index{Finite Euler sums}
\begin{equation} \label{23.56}
\sum_{k=1}^{n} k^3 z^k
\end{equation}\index{Finite Euler sums}
\begin{equation} \nonumber
= \frac{n^3 z^{n+4} +(3n^3+6n^2-4) z^{n+2} -(3n^3+3n^2-3n+1) z^{n+3} -(n+1)^3 z^{n+1} +z^3 + 4z^2 + z}{(1-z)^4},
\end{equation}\index{Finite Euler sums}
\begin{equation} \label{23.57}
\sum_{k=1}^{n} k^4 z^k
= - \frac{n^4z^{n+5}+(-4n^4-4n^3+6n^2-4n+1)z^{n+4} +(6n^4+12n^3-6n^2-12n+11)z^{n+3} }{(1-z)^5}
\end{equation}\index{Finite Euler sums}
\begin{equation} \nonumber
- \frac{(-4n^4-12n^3-6n^2+12n+11)z^{n+2} + (n+1)^4 z^{n+1} - z^4 - 11z^3 - 11z^2 - z}{(1-z)^5}.
\end{equation}\index{Finite Euler sums}

These kind of finite sums lend themselves to substitution in the hyperpyramid identities; in fact for both the square and skewed hyperpyramid versions. It is illustrative here to give a few examples in 2D and 3D of the uses for the Euler Sums\index{Euler, L.}\index{Euler sum} in the context of VPV identities here.

\subsection{The $2D$ square hyperpyramid VPV identity.}

It is known that:

If $|y|, |z| < 1$, with $a+b=1$,
 \begin{equation}   \label{23.57a}
    \prod_{\substack{m,n \geq 1 \\ m \leq n ; \, \gcd(m,n)=1}} \left( \frac{1}{1- y^m z^n} \right)^{\frac{1}{m^a n^b}}
    = \exp\left\{\sum_{n=1}^{\infty} \left( \sum_{m=1}^{n} \frac{y^m}{m^a} \right) \frac{z^n}{n^b} \right\}.
  \end{equation}

Let us set by definition $\gcd(m,n):=(m,n)$ and apply equations (\ref{23.50}) through to (\ref{23.53}) substituted respectively into (\ref{23.57a}). This corresponds to the cases of (\ref{23.57a}) with:

 $y=1$, $a=-1$, implying from condition $a+b=1$ that $b=2$;

 $y=1$, $a=-2$, implying from condition $a+b=1$ that $b=3$;

 $y=1$, $a=-3$, implying from condition $a+b=1$ that $b=4$;

 $y=1$, $a=-4$, implying from condition $a+b=1$ that $b=5$;

 \begin{flushleft}
 so the resulting substitutions give us respectively the four equations
 \end{flushleft}

 \begin{equation}   \nonumber
    \prod_{\substack{m,n \geq 1 \\ m \leq n ; \, (m,n)=1}} \left( \frac{1}{1-z^n} \right)^{\frac{m^1 }{n^2}}
    = \exp\left\{ \sum_{n=1}^{\infty}  \left( \frac{n}{2} + \frac{n^2}{2} \right) \frac{z^n}{n^2} \right\}
\end{equation}
 \begin{equation}   \nonumber
    = \exp\left\{ \frac{1}{2} \log\left(\frac{1}{1-z}\right) + \frac{1}{2}\frac{z}{1-z} \right\},
\end{equation}

 \begin{equation}   \nonumber
    \prod_{\substack{m,n \geq 1 \\ m \leq n ; \, (m,n)=1}} \left( \frac{1}{1-z^n} \right)^{\frac{m^2 }{n^3}}
    = \exp\left\{ \sum_{n=1}^{\infty}  \left( \frac{n}{6} + \frac{n^2}{2} + \frac{n^3}{3} \right) \frac{z^n}{n^3} \right\}
\end{equation}
  \begin{equation}   \nonumber
    = \exp\left\{\frac{1}{6} Li_2(z) + \frac{1}{2} \log\left(\frac{1}{1-z}\right) + \frac{1}{3}\frac{z}{1-z} \right\},
\end{equation}

 \begin{equation}   \nonumber
    \prod_{\substack{m,n \geq 1 \\ m \leq n ; \, (m,n)=1}} \left( \frac{1}{1-z^n} \right)^{\frac{m^3 }{n^4}}
    = \exp\left\{ \sum_{n=1}^{\infty}  \left( \frac{n^2}{4} + \frac{n^3}{3} + \frac{n^4}{4} \right) \frac{z^n}{n^4} \right\}
\end{equation}
 \begin{equation}   \nonumber
    = \exp\left\{\frac{1}{4} Li_2(z) + \frac{1}{3} \log\left(\frac{1}{1-z}\right) + \frac{1}{4}\frac{z}{1-z} \right\},
\end{equation}

 \begin{equation}   \nonumber
    \prod_{\substack{m,n \geq 1 \\ m \leq n ; \, (m,n)=1}} \left( \frac{1}{1-z^n} \right)^{\frac{m^4 }{n^5}}
    = \exp\left\{ \sum_{n=1}^{\infty}  \left( \frac{-n}{30} + \frac{n^3}{3} + \frac{n^4}{2} + \frac{n^5}{5} \right) \frac{z^n}{n^5} \right\}
\end{equation}
 \begin{equation}   \nonumber
    = \exp\left\{\frac{-1}{30} Li_3(z) + \frac{1}{3} \log\left(\frac{1}{1-z}\right) + \frac{1}{2}\frac{z}{(1-z)^2}
    + \frac{1}{5}\frac{z}{1-z} \right\}.
\end{equation}

Based on the above equations and their reciprocal equations, we may assert the following four theorems.

\begin{theorem}    \label{23.57b}
For $|z|<1$,
 \begin{equation}   \nonumber
    \prod_{\substack{m,n \geq 1 \\ m \leq n ; \, (m,n)=1}} \left( \frac{1}{1-z^n} \right)^{\frac{m^1 }{n^2}}
    = \sqrt{\frac{1}{1-z}} \; \exp\left\{ \frac{1}{2}\frac{z}{1-z} \right\},
\end{equation}
 \begin{equation}   \nonumber
    \prod_{\substack{m,n \geq 1 \\ m \leq n ; \, (m,n)=1}} \left(1-z^n\right)^{\frac{m^1 }{n^2}}
    = \sqrt{1-z} \; \exp\left\{ \frac{-1}{2}\frac{z}{1-z} \right\}.
\end{equation}
\end{theorem}

\begin{theorem}   \label{23.57c}
For $|z|<1$,
 \begin{equation}   \nonumber
    \prod_{\substack{m,n \geq 1 \\ m \leq n ; \, (m,n)=1}} \left( \frac{1}{1-z^n} \right)^{\frac{m^2 }{n^3}}
    = \sqrt{\frac{1}{1-z}} \; \exp\left\{\frac{1}{6} Li_2(z) + \frac{1}{3}\frac{z}{1-z} \right\},
\end{equation}
 \begin{equation}   \nonumber
    \prod_{\substack{m,n \geq 1 \\ m \leq n ; \, (m,n)=1}} \left( 1-z^n\right)^{\frac{m^2 }{n^3}}
    = \sqrt{1-z} \; \exp\left\{\frac{-1}{6} Li_2(z) - \frac{1}{3}\frac{z}{1-z} \right\}.
\end{equation}
\end{theorem}

\begin{theorem}   \label{23.57d}
For $|z|<1$,
 \begin{equation}   \nonumber
    \prod_{\substack{m,n \geq 1 \\ m \leq n ; \, (m,n)=1}} \left( \frac{1}{1-z^n} \right)^{\frac{m^3 }{n^4}}
    = \sqrt[3]{\frac{1}{1-z}} \; \exp\left\{\frac{1}{4} Li_2(z) + \frac{1}{4}\frac{z}{1-z} \right\},
\end{equation}
 \begin{equation}   \nonumber
    \prod_{\substack{m,n \geq 1 \\ m \leq n ; \, (m,n)=1}} \left( 1-z^n \right)^{\frac{m^3 }{n^4}}
    = \sqrt[3]{1-z} \; \exp\left\{\frac{-1}{4} Li_2(z) - \frac{1}{4}\frac{z}{1-z} \right\}.
\end{equation}
\end{theorem}

\begin{theorem}   \label{23.57e}
For $|z|<1$,
 \begin{equation}   \nonumber
    \prod_{\substack{m,n \geq 1 \\ m \leq n ; \, (m,n)=1}} \left( \frac{1}{1-z^n} \right)^{\frac{m^4 }{n^5}}
    = \sqrt[3]{\frac{1}{1-z}} \; \exp\left\{\frac{z(7-2z)}{10(1-z)^2} - \frac{1}{30} Li_3(z) \right\},
\end{equation}
 \begin{equation}   \nonumber
    \prod_{\substack{m,n \geq 1 \\ m \leq n ; \, (m,n)=1}} \left( 1-z^n \right)^{\frac{m^4 }{n^5}}
    = \sqrt[3]{1-z} \; \exp\left\{\frac{z(2z-7)}{10(1-z)^2} + \frac{1}{30} Li_3(z) \right\}.
\end{equation}
\end{theorem}

The above four theorems give us a single variable equation in each case that encodes a statement about “weighted integer partitions" of a kind not normally discussed in the partition literature. We could construct a step by step table for the function of the partition, much in the way this was done for the partition congruence function\index{Congruences of partitions}\index{Partition congruences} called the Crank. However equations (\ref{23.54}) through to (\ref{23.57}) substituted respectively into (\ref{23.57a}) can supply us with two-variable 2D generalizations of theorem \ref{23.57b} to theorem \ref{23.57e}. This corresponds to the cases of (\ref{23.57a}) with:

 $a=-1$, implying from condition $a+b=1$ that $b=2$;

 $a=-2$, implying from condition $a+b=1$ that $b=3$;

 $a=-3$, implying from condition $a+b=1$ that $b=4$;

 $a=-4$, implying from condition $a+b=1$ that $b=5$;

 \begin{flushleft}
 so the resulting substitutions (putting $y$ for $z$ in each case) give us respectively the four 2D equation workings,
 \end{flushleft}

 \begin{equation}   \nonumber
    \prod_{\substack{m,n \geq 1 \\ m \leq n ; \, (m,n)=1}} \left( \frac{1}{1- y^m z^n} \right)^{\frac{m^1 }{n^2}}
    = \exp\left\{ \sum_{n=1}^{\infty}  \left( \frac{y - (1 + n) y^{n+1} + n y^{n+2}}{(1-y)^2} \right) \frac{z^n}{n^2} \right\}
\end{equation}
 \begin{equation}   \nonumber
    = \exp\left\{ \frac{y Li_2(z)}{(1-y)^2} - \frac{y Li_2(yz)}{(1-y)^2} - \frac{y}{(1-y)^2} \log\left(\frac{1}{1-yz}\right)
     + \frac{y^2}{(1-y)^2} \log\left(\frac{1}{1-yz}\right) \right\}
\end{equation}
 \begin{equation}   \nonumber
    = \left(1-yz\right)^{\frac{y}{1-y}} \exp\left\{ \frac{y}{(1-y)^2} \left(Li_2(z)- Li_2(yz)\right) \right\};
\end{equation}

 \begin{equation}   \nonumber
    \prod_{\substack{m,n \geq 1 \\ m \leq n ; \, (m,n)=1}} \left( \frac{1}{1- y^m z^n} \right)^{\frac{m^2 }{n^3}}
    \end{equation}
  \begin{equation}   \nonumber
    = \exp\left\{ \sum_{n=1}^{\infty}  \left(\frac{ - n^2 y^{n+3} + (2n^2+2n-1)y^{n+2} -(n^2+2n+1) y^{n+1} + y^2 + y}{(1-y)^3}\right) \frac{z^n}{n^3} \right\}
\end{equation}
\begin{equation}   \nonumber
    = \exp\left\{ \sum_{n=1}^{\infty}  \left(\frac{ -y(y+1) (y^n -1) + 2 n (y-1) y^{n+1} - n^2 (y-1)^2 y^{n+1} }{(1-y)^3}\right) \frac{z^n}{n^3} \right\}
\end{equation}
  \begin{equation}   \nonumber
    = \left(\frac{1}{1-yz}\right)^{\frac{y}{1-y}} \exp\left\{ \frac{y(y+1)}{(1-y)^3}(Li_3(z) -Li_3(yz)) - \frac{2 y}{(1-y)^2} Li_2(yz) \right\};
\end{equation}

 \begin{equation}   \nonumber
    \prod_{\substack{m,n \geq 1 \\ m \leq n ; \, (m,n)=1}} \left( \frac{1}{1- y^m z^n} \right)^{\frac{m^3 }{n^4}}
    = \exp\left\{ \sum_{n=1}^{\infty}  \left( \frac{n^3 y^{n+4} +(3n^3+6n^2-4) y^{n+2}}{(1-y)^4} \right) \frac{z^n}{n^4} \right\}
\end{equation}
 \begin{equation}   \nonumber
    \times \exp\left\{ \sum_{n=1}^{\infty}  \left( \frac{-(n+1)^3 y^{n+1} +y^3 + 4y^2 + y}{(1-y)^4} \right) \frac{z^n}{n^4} \right\}
\end{equation}
 \begin{equation}   \nonumber
    = \exp\left\{ \sum_{n=1}^{\infty}  \left( \frac{-y (y^2 + 4 y + 1) (y^n - 1) + 3 n (y^2 -1) y^{n+1} - 3 n^2 ((y-1)^2 y^{n+1}) + n^3 (y-1)^3 y^{n+1}}{(1-y)^4} \right) \frac{z^n}{n^4} \right\}
\end{equation}
 \begin{equation}   \nonumber
    = \left( 1-yz \right)^{\frac{y}{1-y}} \exp\left\{  \frac{y (y^2 +4y +1)}{(1-y)^4}(Li_4(z)-Li_4(yz)) + \frac{3y(1+y)}{(1-y)^3} Li_3(yz) + \frac{3y}{(1-y)^2}  Li_2(yz) \right\};
\end{equation}

 \begin{equation}   \nonumber
    \prod_{\substack{m,n \geq 1 \\ m \leq n ; \, (m,n)=1}} \left( \frac{1}{1- y^m z^n} \right)^{\frac{m^4 }{n^5}}
\end{equation}
 \begin{equation}   \nonumber
    = \exp\left\{ \sum_{n=1}^{\infty}  \left( - \frac{n^4y^{n+5}+(-4n^4-4n^3+6n^2-4n+1)y^{n+4} }{(1-y)^5} \right) \frac{z^n}{n^5} \right\}
\end{equation}
 \begin{equation}   \nonumber
    \times \exp\left\{ \sum_{n=1}^{\infty}  \left( - \frac{(6n^4+12n^3-6n^2-12n+11)y^{n+3} }{(1-y)^5} \right) \frac{z^n}{n^5} \right\}
\end{equation}
 \begin{equation}   \nonumber
    \times \exp\left\{ \sum_{n=1}^{\infty}  \left( - \frac{(-4n^4-12n^3-6n^2+12n+11)y^{n+2} + (n+1)^4 y^{n+1} }{(1-y)^5} \right) \frac{z^n}{n^5} \right\}
\end{equation}
 \begin{equation}   \nonumber
    \times \exp\left\{ \sum_{n=1}^{\infty}  \left( + \frac{y^4 + 11y^3 + 11y^2 + z}{(1-y)^5} \right) \frac{z^n}{n^5} \right\}
\end{equation}

\begin{equation}   \nonumber
    = \exp\left\{ \frac{y (y + 1) (y^2 + 10 y + 1)}{(1-y)^5}(Li_5(yz) - Li_5(z)) \right\}
\end{equation}
\begin{equation}   \nonumber
    \times \exp\left\{ \frac{4}{(1-y)^5} ((-y^4 +3y^3 +y)Li_4(yz) + 3Li_4(z)) \right\}
\end{equation}
 \begin{equation}   \nonumber
    \times \exp\left\{ \frac{6}{(1-y)^5} ((y^4-y^3+y) Li_3(yz) - Li_3(z))\right\}
\end{equation}
\begin{equation}   \nonumber
    \times \exp\left\{\frac{4}{(1-y)^5}((y-3y^3-y^4) Li_2(yz) - 3Li_2(z)) \right\}
\end{equation}
\begin{equation}   \nonumber
    \times \exp\left\{ \frac{1}{(1-y)^5}((y^5 - 4 y^4 + 6y^3 +y) Li_1(yz) - 4Li_1(z))\right\}
\end{equation}

 \begin{equation}   \nonumber
    = \frac{\left(1-z\right)^{\frac{4}{(1-y)^5}}}{\left(1-yz\right)^{\frac{y^5 - 4 y^4 + 6y^3 +y}{(1-y)^5}}}
    \end{equation}
\begin{equation}   \nonumber
    \times \exp\left\{ \frac{y (y + 1) (y^2 + 10 y + 1)}{(1-y)^5}(Li_5(yz) - Li_5(z)) \right\}
\end{equation}
\begin{equation}   \nonumber
    \times \exp\left\{ \frac{4}{(1-y)^5} ((-y^4 +3y^3 +y)Li_4(yz) + 3Li_4(z)) \right\}
\end{equation}
 \begin{equation}   \nonumber
    \times \exp\left\{ \frac{6}{(1-y)^5} ((y^4-y^3+y) Li_3(yz) - Li_3(z))\right\}
\end{equation}
\begin{equation}   \nonumber
    \times \exp\left\{\frac{4}{(1-y)^5}((y-3y^3-y^4) Li_2(yz) - 3Li_2(z)) \right\}.
\end{equation}

\bigskip

Based on the above equations and their reciprocal equations, we may assert the following four 2D theorems.

\begin{theorem}    \label{23.57f}
For $|z|<1$,
 \begin{equation}   \nonumber
    \prod_{\substack{m,n \geq 1 \\ m \leq n ; \, (m,n)=1}} \left( \frac{1}{1- y^m z^n} \right)^{\frac{m^1 }{n^2}}
    = \left(1-yz\right)^{\frac{y}{1-y}} \exp\left\{ \frac{y}{(1-y)^2} \left(Li_2(z)- Li_2(yz)\right) \right\},
\end{equation}
 \begin{equation}   \nonumber
    \prod_{\substack{m,n \geq 1 \\ m \leq n ; \, (m,n)=1}} \left(1- y^m z^n\right)^{\frac{m^1 }{n^2}}
    = \left(\frac{1}{1-yz}\right)^{\frac{y}{1-y}} \exp\left\{- \frac{y}{(1-y)^2} \left(Li_2(z)- Li_2(yz)\right) \right\}.
\end{equation}
\end{theorem}

\begin{theorem}   \label{23.57g}
For $|z|<1$,
 \begin{equation}   \nonumber
    \prod_{\substack{m,n \geq 1 \\ m \leq n ; \, (m,n)=1}} \left( \frac{1}{1- y^m z^n} \right)^{\frac{m^2 }{n^3}}
    = \left(\frac{1}{1-yz}\right)^{\frac{y}{1-y}}
    \end{equation}
 \begin{equation}   \nonumber
    \times \exp\left\{ \frac{y(y+1)}{(1-y)^3}(Li_3(z) -Li_3(yz)) - \frac{2 y}{(1-y)^2} Li_2(yz) \right\},
\end{equation}
 \begin{equation}   \nonumber
    \prod_{\substack{m,n \geq 1 \\ m \leq n ; \, (m,n)=1}} \left( 1- y^m z^n\right)^{\frac{m^2 }{n^3}}
    = \left(1-yz\right)^{\frac{y}{1-y}}
    \end{equation}
 \begin{equation}   \nonumber
    \times \exp\left\{- \frac{y(y+1)}{(1-y)^3}(Li_3(z) -Li_3(yz)) + \frac{2 y}{(1-y)^2} Li_2(yz) \right\},
\end{equation}
\end{theorem}

\begin{theorem}   \label{23.57h}
For $|z|<1$,
 \begin{equation}   \nonumber
    \prod_{\substack{m,n \geq 1 \\ m \leq n ; \, (m,n)=1}} \left( \frac{1}{1- y^m z^n} \right)^{\frac{m^3 }{n^4}}
    = \left( 1-yz \right)^{\frac{y}{1-y}}
 \end{equation}
 \begin{equation}   \nonumber
     \times \exp\left\{  \frac{y (y^2 +4y +1)}{(1-y)^4}(Li_4(z)-Li_4(yz)) + \frac{3y(1+y)}{(1-y)^3} Li_3(yz) + \frac{3y}{(1-y)^2}  Li_2(yz) \right\},
\end{equation}
 \begin{equation}   \nonumber
    \prod_{\substack{m,n \geq 1 \\ m \leq n ; \, (m,n)=1}} \left( 1- y^m z^n \right)^{\frac{m^3 }{n^4}}
    = \left(\frac{1}{1-yz} \right)^{\frac{y}{1-y}}
    \end{equation}
 \begin{equation}   \nonumber
 \times \exp\left\{ -\frac{y (y^2 +4y +1)}{(1-y)^4}(Li_4(z)-Li_4(yz)) - \frac{3y(1+y)}{(1-y)^3} Li_3(yz) - \frac{3y}{(1-y)^2}  Li_2(yz) \right\}.
\end{equation}
\end{theorem}

\begin{theorem}   \label{23.57i}
For $|z|<1$,
 \begin{equation}   \nonumber
    \prod_{\substack{m,n \geq 1 \\ m \leq n ; \, (m,n)=1}} \left( \frac{1}{1- y^m z^n} \right)^{\frac{m^4 }{n^5}}
    = \left\{\frac{\left(1-z\right)^{4}}{\left(1-yz\right)^{y^5 - 4 y^4 + 6y^3 +y}}\right\}^{\frac{1}{(1-y)^5}}
    \end{equation}
\begin{equation}   \nonumber
    \times \exp\left\{ \frac{y (y + 1) (y^2 + 10 y + 1)}{(1-y)^5}(Li_5(yz) - Li_5(z)) \right\}
\end{equation}
\begin{equation}   \nonumber
    \times \exp\left\{ \frac{4}{(1-y)^5} ((-y^4 +3y^3 +y)Li_4(yz) + 3Li_4(z)) \right\}
\end{equation}
 \begin{equation}   \nonumber
    \times \exp\left\{ \frac{6}{(1-y)^5} ((y^4-y^3+y) Li_3(yz) - Li_3(z))\right\}
\end{equation}
\begin{equation}   \nonumber
    \times \exp\left\{\frac{4}{(1-y)^5}((y-3y^3-y^4) Li_2(yz) - 3Li_2(z)) \right\},
\end{equation}

\begin{equation}   \nonumber
    \prod_{\substack{m,n \geq 1 \\ m \leq n ; \, (m,n)=1}} \left( 1- y^m z^n \right)^{\frac{m^4 }{n^5}}
    = \left\{\frac{\left(1-yz\right)^{y^5 - 4 y^4 + 6y^3 +y}}{\left(1-z\right)^{4}}\right\}^{\frac{1}{(1-y)^5}}
    \end{equation}
\begin{equation}   \nonumber
    \times \exp\left\{ \frac{-y (y + 1) (y^2 + 10 y + 1)}{(1-y)^5}(Li_5(yz) - Li_5(z)) \right\}
\end{equation}
\begin{equation}   \nonumber
    \times \exp\left\{ \frac{-4}{(1-y)^5} ((-y^4 +3y^3 +y)Li_4(yz) + 3Li_4(z)) \right\}
\end{equation}
 \begin{equation}   \nonumber
    \times \exp\left\{ \frac{-6}{(1-y)^5} ((y^4-y^3+y) Li_3(yz) - Li_3(z))\right\}
\end{equation}
\begin{equation}   \nonumber
    \times \exp\left\{\frac{-4}{(1-y)^5}((y-3y^3-y^4) Li_2(yz) - 3Li_2(z)) \right\}.
\end{equation}
\end{theorem}

\subsection{The $3D$ square hyperpyramid VPV identity.}

It is known that:

If $|x|, |y|, |z| < 1$, with $a+b+c=1$,
 \begin{equation}   \label{23.58}
    \prod_{\substack{l,m,n \geq 1 \\ l,m \leq n ; \, \gcd(l,m,n)=1}} \left( \frac{1}{1-x^l y^m z^n} \right)^{\frac{1}{l^a m^b n^c}}
    = \exp\left\{ \sum_{n=1}^{\infty}  \left( \sum_{l=1}^{n} \frac{x^l}{l^a} \right) \left( \sum_{m=1}^{n} \frac{y^m}{m^b} \right) \frac{z^n}{n^c} \right\}.
  \end{equation}

In one of the simpler examples we can give, let us apply equations (\ref{23.50}) and (\ref{23.51}) to (\ref{23.58}). This corresponds to the case of (\ref{23.58}) with $x=y=1$, $a=-1$, $b=-2$, implying from condition $a+b+c=1$ that $c=4$. The resulting substitution is

 \begin{equation}   \nonumber
    \prod_{\substack{l,m,n \geq 1 \\ l,m \leq n ; \, \gcd(l,m,n)=1}} \left( \frac{1}{1-z^n} \right)^{\frac{ l m^2 }{n^4}}
    = \exp\left\{ \sum_{n=1}^{\infty}  \left( \frac{n}{2} + \frac{n^2}{2} \right) \left( \frac{n}{6} + \frac{n^2}{2} + \frac{n^3}{3} \right) \frac{z^n}{n^4} \right\}
\end{equation}
\begin{equation} \nonumber
= \exp\left\{ \sum_{n=1}^{\infty} \left( \frac{n^2}{12} + \frac{n^3}{3} + \frac{5 n^4}{12} + \frac{n^5}{6} \right) \frac{z^n}{n^4} \right\}
\end{equation}
\begin{equation} \nonumber
= \exp\left\{ \frac{1}{12}Li_2(z) + \frac{1}{3}Li_1(z) + \frac{5}{12}Li_0(z) + \frac{1}{6}Li_{-1}(z)  \right\}
\end{equation}
\begin{equation} \nonumber
= \exp\left\{ \frac{1}{12}Li_2(z) + \frac{1}{3} \log\left( \frac{1}{1-z} \right) + \frac{5}{12} \frac{z}{1-z} + \frac{1}{6} \frac{z}{(1-z)^2}  \right\}
\end{equation}
\begin{equation} \nonumber
= \exp\left\{ \frac{1}{12}Li_2(z) + \frac{z(7-5z))}{12(1-z)^2} + \frac{1}{3} \log\left( \frac{1}{1-z} \right) \right\}
\end{equation}
\begin{equation} \nonumber
=    \sqrt[3]{\left( \frac{1}{1-z} \right)} \;  \exp\left\{ \frac{1}{12}Li_2(z) + \frac{z(7-5z))}{12(1-z)^2} \right\}.
\end{equation}

Therefore we have proven that if $|z|<1$,

 \begin{equation}   \label{23.59}
    \prod_{\substack{l,m,n \geq 1 \\ l,m \leq n ; \, \gcd(l,m,n)=1}} \left( \frac{1}{1-z^n} \right)^{\frac{ l m^2 }{n^4}}
    =    \sqrt[3]{\left( \frac{1}{1-z} \right)} \;  \exp\left\{ \frac{1}{12} \left( Li_2(z) + \frac{z(7-5z)}{(1-z)^2}\right) \right\},
\end{equation}

and the equivalent result,

 \begin{equation}   \label{23.60}
    \prod_{\substack{l,m,n \geq 1 \\ l,m \leq n ; \, \gcd(l,m,n)=1}} \left( 1-z^n \right)^{\frac{ l m^2 }{n^4}}
    =    \sqrt[3]{\left( 1-z \right)} \;  \exp\left\{ - \frac{1}{12} \left( Li_2(z) + \frac{z(7-5z)}{(1-z)^2}\right) \right\}.
\end{equation}

There is a corresponding example we can give, if we apply equations (\ref{23.54}) and (\ref{23.55}) to (\ref{23.58}). This corresponds to the case of (\ref{23.58}) with $a=-1$, $b=-2$, implying from condition $a+b+c=1$ that $c=4$. The resulting substitution gives us

 \begin{equation}     \label{23.61}
    \prod_{\substack{l,m,n \geq 1 \\ l,m \leq n ; \, \gcd(l,m,n)=1}} \left( \frac{1}{1-x^l y^m z^n} \right)^{\frac{ l m^2 }{n^4}}
\end{equation}
\begin{equation} \nonumber
    = \exp\left\{ \sum_{n=1}^{\infty} A(n,x) B(n,y) \frac{z^n}{n^4} \right\}
\end{equation}
where
\begin{equation} \nonumber
 A(n,x)   = \frac{x - (1 + n) x^{n+1} + n x^{n+2}}{(1-x)^2}
\end{equation}
and
\begin{equation} \nonumber
 B(n,y)   =  \frac{ - n^2 y^{n+3} + (2n^2+2n-1)y^{n+2} -(n^2+2n+1) y^{n+1} + y^2 + y}{(1-y)^3},
\end{equation}

which, after expanding and simplifying is equivalent to the 

\begin{theorem} \label{23.61a}

 \begin{equation}     \label{23.62}
    \prod_{\substack{l,m,n \geq 1 \\ l,m \leq n ; \, \gcd(l,m,n)=1}} \left( \frac{1}{1-x^l y^m z^n} \right)^{\frac{ l m^2 }{n^4}}
\end{equation}

\begin{equation} \nonumber
    =  \left(\frac{1}{1-xyz} \right)^{\frac{xy}{(1-x)^2(1-y)^3}}    \times \exp\left\{ \frac{-xy(2x+y-3)}{(1-x)^3(1-y)^4} Li_2(xyz) \right\}
\end{equation}
\begin{equation} \nonumber
    \times \exp\left\{  \frac{xy}{(1-x)^3(1-y)^3} Li_2(yz) + \frac{xy(y+1)}{(1-x)^3(1-y)^5} Li_4(xyz) \right\}
\end{equation}
\begin{equation} \nonumber
    \times \exp\left\{ \frac{-xy(xy+x+y-3)}{(1-x)^3(1-y)^5} Li_3(xyz) - \frac{xy(y+1)}{(1-x)^2(1-y)^5}Li_3(xz) \right\}
\end{equation}
\begin{equation} \nonumber
    \times \exp\left\{ - \frac{xy(y+1)}{(1-x)^3(1-y)^5} Li_4(x z) - \frac{2xy}{(1-x)^3(1-y)^4}Li_3(yz)  \right\}.
\end{equation}
\begin{equation} \nonumber
    \times \exp\left\{ \frac{-xy(1+y)}{(1-x)^3(1-y)^5}Li_4(yz) + \frac{xy(1+y) }{(1-x)^3(1-y)^5}Li_4(z) \right\}.
\end{equation}
\end{theorem}

\textbf{Proof of theorem \ref{23.61a}:}
The necessary workings and simplification in going from (\ref{23.61}) to arrive at (\ref{23.62}) are given in the following page and a half of analysis.

 \begin{equation}   \nonumber
    \prod_{\substack{l,m,n \geq 1 \\ l,m \leq n ; \, \gcd(l,m,n)=1}} \left( \frac{1}{1-x^l y^m z^n} \right)^{\frac{ l m^2 }{n^4}}
\end{equation}

\begin{equation} \nonumber
    = \exp\left\{ \frac{1}{(1-x)^3(1-y)^5} \sum_{n=1}^{\infty} (n^3x^{n+1}y^{n+1}-2n^3x^{n+1}y^{n+2}+n^3x^{n+1}y^{n+3}) \frac{z^n}{n^4} \right\}
\end{equation}
\begin{equation} \nonumber
    \times \exp\left\{ \frac{1}{(1-x)^3(1-y)^5} \sum_{n=1}^{\infty} (-n^3x^{n+2}y^{n+1}+2n^3x^{n+2}y^{n+2}-n^3x^{n+2}y^{n+3}) \frac{z^n}{n^4} \right\}
\end{equation}
\begin{equation} \nonumber
    \times \exp\left\{ \frac{1}{(1-x)^3(1-y)^5} \sum_{n=1}^{\infty} (+3n^2x^{n+1}y^{n+1}-4n^2x^{n+1}y^{n+2}+n^2x^{n+1}y^{n+3}) \frac{z^n}{n^4} \right\}
\end{equation}
\begin{equation} \nonumber
    \times \exp\left\{ \frac{1}{(1-x)^3(1-y)^5} \sum_{n=1}^{\infty} (-2n^2x^{n+2}y^{n+1}+2n^2x^{n+2}y^{n+2}-n^2xy^{n+1}) \frac{z^n}{n^4} \right\}
\end{equation}
\begin{equation} \nonumber
    \times \exp\left\{ \frac{1}{(1-x)^3(1-y)^5} \sum_{n=1}^{\infty} (+2n^2xy^{n+2}-n^2xy^{n+3}+3nx^{n+1}y^{n+1}) \frac{z^n}{n^4} \right\}
\end{equation}
\begin{equation} \nonumber
    \times \exp\left\{ \frac{1}{(1-x)^3(1-y)^5} \sum_{n=1}^{\infty} (+x^{n+1}y^{n+1}-nx^{n+1}y^{n+2}+x^{n+1}y^{n+2}) \frac{z^n}{n^4} \right\}
\end{equation}
\begin{equation} \nonumber
    \times \exp\left\{ \frac{1}{(1-x)^3(1-y)^5} \sum_{n=1}^{\infty} (-nx^{n+2}y^{n+1}-nx^{n+2}y^{n+2}-ny^2x^{n+1}-y^2x^{n+1}) \frac{z^n}{n^4} \right\}
\end{equation}
\begin{equation} \nonumber
    \times \exp\left\{ \frac{1}{(1-x)^3(1-y)^5} \sum_{n=1}^{\infty} (+ny^2x^{n+2}-nyx^{n+1}-yx^{n+1}+nyx^{n+2}-2nxy^{n+1}) \frac{z^n}{n^4} \right\}
\end{equation}
\begin{equation} \nonumber
    \times \exp\left\{ \frac{1}{(1-x)^3(1-y)^5} \sum_{n=1}^{\infty} (-xy^{n+1}+2nxy^{n+2}-xy^{n+2}+xy^2+xy) \frac{z^n}{n^4} \right\}
\end{equation}

which works out through the following analysis to

 \begin{equation}   \nonumber
    \prod_{\substack{l,m,n \geq 1 \\ l,m \leq n ; \, \gcd(l,m,n)=1}} \left( \frac{1}{1-x^l y^m z^n} \right)^{\frac{ l m^2 }{n^4}}
\end{equation}

\begin{equation} \nonumber
    = \exp\left\{ \frac{1}{(1-x)^3(1-y)^5}  (xyLi_1(xyz)-2xy^2Li_1(xyz)+xy^3Li_1(xyz))  \right\}
\end{equation}
\begin{equation} \nonumber
    \times \exp\left\{ \frac{1}{(1-x)^3(1-y)^5} (-x^2yLi_1(xyz)+2x^2y^2Li_1(xyz)-x^2y^3Li_1(xyz)) \right\}
\end{equation}
\begin{equation} \nonumber
    \times \exp\left\{ \frac{1}{(1-x)^3(1-y)^5} (+3xyLi_2(xyz)-4xy^2Li_2(xyz)+xy^3Li_2(xyz)) \right\}
\end{equation}
\begin{equation} \nonumber
    \times \exp\left\{ \frac{1}{(1-x)^3(1-y)^5} (-2x^2yLi_2(xyz)+2x^2y^2Li_2(xyz)-xyLi_2(yz)) \right\}
\end{equation}
\begin{equation} \nonumber
    \times \exp\left\{ \frac{1}{(1-x)^3(1-y)^5} (+2xy^2Li_2(yz)-xy^3Li_2(yz)+3xyLi_3(xyz)) \right\}
\end{equation}
\begin{equation} \nonumber
    \times \exp\left\{ \frac{1}{(1-x)^3(1-y)^5} (+xyLi_4(xyz)-xy^2Li_3(xyz)+xy^2Li_4(xyz)) \right\}
\end{equation}
\begin{equation} \nonumber
    \times \exp\left\{ \frac{1}{(1-x)^3(1-y)^5} (-x^2yLi_3(xyz)-x^2y^2Li_3(xyz)-xy^2Li_3(xz)-xy^2Li_4(xz)) \right\}
\end{equation}
\begin{equation} \nonumber
    \times \exp\left\{ \frac{1}{(1-x)^3(1-y)^5} (+x^2y^2Li_3(xz)-xyLi_3(xz)-xyLi_4(xz)+x^2yLi_3(xz)) \right\}
\end{equation}
\begin{equation} \nonumber
    \times \exp\left\{ \frac{1}{(1-x)^3(1-y)^5} (-2xyLi_3(yz) -xyLi_4(yz)+2xy^2Li_3(yz)-xy^2Li_4(yz)) \right\}
\end{equation}
\begin{equation} \nonumber
    \times \exp\left\{ \frac{1}{(1-x)^3(1-y)^5} (+xy^2Li_4(z)+xyLi_4(z)) \right\}
\end{equation}

\begin{equation} \nonumber
    = \exp\left\{ \frac{1}{(1-x)^3(1-y)^5}  ((x-1)x(y - 1)^2y \log(1-xyz))  \right\}
\end{equation}
\begin{equation} \nonumber
    \times \exp\left\{ \frac{1}{(1-x)^3(1-y)^5} (x(y-1)y(2x+y-3) Li_2(xyz)) \right\}
\end{equation}
\begin{equation} \nonumber
    \times \exp\left\{ \frac{1}{(1-x)^3(1-y)^5} (-x(y-1)^2 y Li_2(yz)) \right\}
\end{equation}
\begin{equation} \nonumber
    \times \exp\left\{ \frac{1}{(1-x)^3(1-y)^5} (xy(y + 1) Li_4(xyz)) \right\}
\end{equation}
\begin{equation} \nonumber
    \times \exp\left\{ \frac{1}{(1-x)^3(1-y)^5} (-xy(xy+x+y-3) Li_3(xyz)) \right\}
\end{equation}
\begin{equation} \nonumber
    \times \exp\left\{ \frac{1}{(1-x)^3(1-y)^5} ((x-1)xy(y+1)Li_3(xz)-xy(y+1)Li_4(x z)) \right\}
\end{equation}
\begin{equation} \nonumber
    \times \exp\left\{ \frac{1}{(1-x)^3(1-y)^5} (+2x(y-1)yLi_3(y z)-xy(1+y)Li_4(yz)+xy(1+y)Li_4(z)) \right\}
\end{equation}

\begin{equation} \nonumber
    =  \left(\frac{1}{1-xyz} \right)^{\frac{xy}{(1-x)^2(1-y)^3}}    \times \exp\left\{ \frac{-xy(2x+y-3)}{(1-x)^3(1-y)^4} Li_2(xyz) \right\}
\end{equation}
\begin{equation} \nonumber
    \times \exp\left\{  \frac{xy}{(1-x)^3(1-y)^3} Li_2(yz) + \frac{xy(y+1)}{(1-x)^3(1-y)^5} Li_4(xyz) \right\}
\end{equation}
\begin{equation} \nonumber
    \times \exp\left\{ \frac{-xy(xy+x+y-3)}{(1-x)^3(1-y)^5} Li_3(xyz) - \frac{xy(y+1)}{(1-x)^2(1-y)^5}Li_3(xz) \right\}
\end{equation}
\begin{equation} \nonumber
    \times \exp\left\{ - \frac{xy(y+1)}{(1-x)^3(1-y)^5} Li_4(x z) - \frac{2xy}{(1-x)^3(1-y)^4}Li_3(yz)  \right\}.
\end{equation}
\begin{equation} \nonumber
    \times \exp\left\{ \frac{-xy(1+y)}{(1-x)^3(1-y)^5}Li_4(yz) + \frac{xy(1+y) }{(1-x)^3(1-y)^5}Li_4(z) \right\}.
\end{equation}

$\quad \quad \quad \blacksquare$

The casual reader will see that the earlier simple summations arising from substituted equations (\ref{23.50}) through to (\ref{23.53}) yield much more immediate and simple results than the more elaborate substituted equations (\ref{23.54}) through to (\ref{23.57}).

\section{Parametric Euler Sum Identities}\index{Parametric Euler sum}\index{Euler, L.}

We begin by defining various special functions which we will later use. For
$\Re(x) > 0$, the gamma function and its logarithmic derivative (digamma function) are defined by
\begin{equation} \nonumber
  \Gamma(x) := \int_{0}^{\infty} e^{-t} t^{x-1} dt, \quad  \Psi(x) := \frac{\Gamma'(x)}{\Gamma(x)}.
\end{equation}

We note that
\begin{equation} \nonumber
  \sum_{n=1}^{\infty} \frac{x}{n(n-x)} := -\Psi(1-x) - \gamma, \quad \textmd{where} \quad \gamma := \lim_{N\rightarrow\infty} \left(\sum_{n=1}^{N}\frac{1}{n} - \log N  \right)
\end{equation}

is Euler's constant\index{Euler, L.}. In addition to these we use the classical Riemann zeta function\index{Riemann zeta function}
\begin{equation} \nonumber
  \zeta(s) := \sum_{n=1}^{\infty} \frac{1}{n^s}  \quad  \Re(s) > 1,
\end{equation}

and extend this to
\begin{equation} \label{23.63}
  \zeta(s,t) := \sum_{n>m>0} \frac{1}{n^s m^t} = \sum_{n=1}^{\infty} \frac{1}{n^s} \sum_{m=1}^{n-1} \frac{1}{m^t} \quad  \Re(s) > 1, \; \Re(t) > 1.
\end{equation}
defining a double Euler sum\index{Parametric Euler sum}\index{Euler, L.}\index{Euler sum}. The two-place function (\ref{23.63}) was first introduced by Euler\index{Euler, L.},
who noted the reflection formula
\begin{equation} \label{23.64}
  \zeta(s,t) + \zeta(t,s) = \zeta(s) \zeta(t) - \zeta(s+t) \quad  \Re(s) > 1, \; \Re(t) > 1,
\end{equation}
and the reduction formula\index{Parametric Euler sum}
\begin{equation} \label{23.65}
  \zeta(s,1) = \frac{1}{2}s \, \zeta(s+1) - \frac{1}{2} \sum_{k=1}^{s-2}\zeta(k+1) \zeta(s-k), \quad \quad 1 < s \in \mathbb{Z} ,
\end{equation}
which reduces the double Euler sum\index{Parametric Euler sum}\index{Euler, L.} $\zeta(s, 1)$ to a sum of products of classical Riemann
$\zeta$-values. From setting $s=2$ in (\ref{23.65}) yields the important sum,
\begin{equation} \label{23.66}
  \zeta(2,1) = \zeta(3).
\end{equation}
The next cases are,
\begin{equation} \label{23.66a}
  \zeta(3,1) = \frac{3}{2} \zeta(4) - \frac{1}{2} \zeta(2),
\end{equation}
\begin{equation} \label{23.66b}
  \zeta(4,1) = 2 \zeta(5) - \frac{1}{2} \zeta(2)^2,
\end{equation}
\begin{equation} \label{23.66c}
  \zeta(5,1) = \frac{5}{2} \zeta(6) - \zeta(2)\zeta(3).
\end{equation}

In 1994 Bailey et al. \cite{dB1994} formulated the following definitions for Euler sums\index{Parametric Euler sum}\index{Euler, L.}:

\begin{eqnarray}
 \label{23.67}  s_h(m,n) &=& \sum_{k=1}^{\infty} \left( \frac{1}{1} + \frac{1}{2} + \cdots + \frac{1}{k} \right)^m \frac{1}{(k+1)^n}, \\
 \label{23.68}  s_a(m,n) &=& \sum_{k=1}^{\infty} \left( \frac{1}{1} - \frac{1}{2} + \cdots - \frac{(-1)^{k+1}}{k} \right)^m \frac{1}{(k+1)^n}, \\
 \label{23.69}  a_h(m,n) &=& \sum_{k=1}^{\infty} \left( \frac{1}{1} + \frac{1}{2} + \cdots + \frac{1}{k} \right)^m\frac{(-1)^{k+1}}{(k+1)^n}, \\
 \label{23.70}  a_a(m,n) &=& \sum_{k=1}^{\infty} \left( \frac{1}{1} - \frac{1}{2} + \cdots + \frac{(-1)^{k+1}}{k} \right)^m \frac{(-1)^{k+1}}{(k+1)^n}, \\
 \label{23.71}  \sigma_h(m,n) &=& \sum_{k=1}^{\infty} \left( \frac{1}{1^m} + \frac{1}{2^m} + \cdots + \frac{1}{k^m} \right) \frac{1}{(k+1)^n}, \\
 \label{23.72}  \sigma_a(m,n) &=& \sum_{k=1}^{\infty} \left( \frac{1}{1^m} - \frac{1}{2^m} + \cdots - \frac{(-1)^{k+1}}{k^m} \right) \frac{1}{(k+1)^n}, \\
 \label{23.73}  \alpha_h(m,n) &=& \sum_{k=1}^{\infty} \left( \frac{1}{1^m} + \frac{1}{2^m} + \cdots + \frac{1}{k^m} \right) \frac{(-1)^{k+1}}{(k+1)^n}, \\
 \label{23.74}  \alpha_a(m,n) &=& \sum_{k=1}^{\infty} \left( \frac{1}{1^m} - \frac{1}{2^m} + \cdots + \frac{(-1)^{k+1}}{k^m} \right) \frac{(-1)^{k+1}}{(k+1)^n},
\end{eqnarray}

and from these come three notable special forms\index{Parametric Euler sum}\index{Euler, L.},
\begin{equation} \label{23.75}
  s_h(m,n) = \sum_{k=1}^{\infty} \left( \gamma + \Psi(k+1) \right)^m \frac{(-1)^{k+1}}{(k+1)^n},
\end{equation}
\begin{equation} \label{23.76}
  a_a(m,n) = \sum_{k=1}^{\infty} \left( \log2
 + \frac{1}{2} (-1)^{k+1}\left[ \Psi(\frac{k}{2}+\frac{1}{2})-\Psi(\frac{k}{2}+1)\right]\right)^m \frac{(-1)^{k+1}}{(k+1)^n},
\end{equation}
\begin{equation} \label{23.77}
  \alpha_h(m,n) = \left(1- \frac{2}{2^{m+n}}\right) \sum_{k=1}^{\infty} \frac{(-1)^{k+1} H_{k,m}}{k^n},
\end{equation}
where $H_{k,m}= \sum_{j=1}^{k} j^{-m}$ is a generalized harmonic number\index{Harmonic number}\index{Generalized harmonic number} and $\Psi(n)= \Gamma'(z)/\Gamma(z)$ is the digamma function\index{Digamma function}.

At this point we note that each Euler sum\index{Parametric Euler sum}\index{Euler, L.}\index{Euler sum} defined by equations (\ref{23.67}) to (\ref{23.74}) and therefore also for (\ref{23.75}) to (\ref{23.77}) when viewed as lattice point multi-dimensional sums have the following characteristics. The first four sums (\ref{23.67}) to (\ref{23.70}) are $m+1$-dimensional lattice sums with $m$ dimensions identical copies of each other, and the other four sums (\ref{23.71}) to (\ref{23.74}) are each 2D lattice sums. Of these, all eight sums are summations taken over radial from origin lattice point regions applicable to our higher dimensional VPV identities. In fact, the above definitions can be extended to power series that also conform to the radial from origin region requirement for applying more general hyperquadrant and hyperpyramid VPV theorems.

There are many well known particular cases of the Euler sums\index{Euler, L.}\index{Euler sum} as given in (\ref{23.67}) to (\ref{23.74}). A few of them for $s_h(m,n)$ were found by experimental means by Bailey et al. \cite{dB1994} and are:
\begin{equation} \label{23.78}
  s_h(3,2) = \frac{15}{2} \zeta(5) + \zeta(2) \zeta(3),
\end{equation}
\begin{equation} \label{23.79}
  s_h(4,3) = -\frac{109}{8} \zeta(7) + \frac{37}{2} \zeta(3) \zeta(4) -5 \zeta(2) \zeta(5),
\end{equation}
\begin{equation} \label{23.80}
  s_h(5,4) = \frac{890}{9}\zeta(9) +66\zeta(4)\zeta(5) -\frac{4295}{24}\zeta(2)\zeta(5) -5\zeta(3)^3 +\frac{265}{8} \zeta(2)\zeta(7),
\end{equation}
each of which are cases of
\begin{equation} \label{23.81}
  s_h(m+1,m) = \sum_{k=1}^{\infty} \left( \frac{1}{1} + \frac{1}{2} + \cdots + \frac{1}{k} \right)^{m+1} \frac{1}{(k+1)^m}.
\end{equation}

The sums (\ref{23.78}) to (\ref{23.80}) are likely to be true as they were verified to several hundred decimal places. All that we need with respect to them is an actual proof.

Equation (\ref{23.81}) is directly applicable to the VPV hyperpyramid identity case, which may be established from methods in Campbell \cite{gC2019}, and is here stated as,
\begin{equation} \label{23.82}
\prod_{\substack{\gcd(j_1,j_2,\ldots,j_{m+1},k)=1 \\ j_1,j_2,\ldots,j_{m+1}<k \\ j_1,j_2,\ldots,j_{m+1} \geq 0; k \geq 1}} \left(\frac{1}{1-y^{j_1+j_2+...+j_{m+1}}z^k} \right)^{\frac{1}{k}} \quad \quad
\end{equation}
\begin{equation*}
  \quad \quad  = \exp\left[ \sum_{k=1}^{\infty} \left( \frac{y}{1} + \frac{y^2}{2} + \cdots + \frac{y^k}{k} \right)^{m+1} \frac{z^{k+1}}{(k+1)^m}\right].
\end{equation*}

Applying our three experimentally obtained sums (\ref{23.78}) to (\ref{23.80}) successively, setting $y\rightarrow1$ and $z\rightarrow1$ in the three identities,
\begin{equation} \label{23.83}
\prod_{\substack{\gcd(j_1,j_2,j_3,k)=1 \\ j_1,j_2,j_3<k \\ j_1,j_2,j_3 \geq 1; k \geq 2}} \left(\frac{1}{1-y^{j_1+j_2+j_3}z^k} \right)^{\frac{1}{k}}
   = \exp\left[ \sum_{k=1}^{\infty} \left( \frac{y}{1} + \frac{y^2}{2} + \cdots + \frac{y^k}{k} \right)^{3} \frac{z^{k+1}}{(k+1)^2}\right].
\end{equation}
\begin{equation} \label{23.84}
\prod_{\substack{\gcd(j_1,j_2,j_3,j_4,k)=1 \\ j_1,j_2,j_3,j_4<k \\ j_1,j_2,j_3,j_4\geq 1; k \geq 2}} \left(\frac{1}{1-y^{j_1+j_2+j_3+j_4}z^k} \right)^{\frac{1}{k}}
   = \exp\left[ \sum_{k=1}^{\infty} \left( \frac{y}{1} + \frac{y^2}{2} + \cdots + \frac{y^k}{k} \right)^{4} \frac{z^{k+1}}{(k+1)^3}\right].
\end{equation}
\begin{equation} \label{23.85}
\prod_{\substack{\gcd(j_1,j_2,j_3,j_4,j_5,k)=1 \\ j_1,j_2,j_3,j_4,j_5<k \\ j_1,j_2,j_3,j_4,j_5\geq 1; k \geq 2}} \left(\frac{1}{1-y^{j_1+j_2+j_3+j_4+j_5}z^k} \right)^{\frac{1}{k}}
   = \exp\left[ \sum_{k=1}^{\infty} \left( \frac{y}{1} + \frac{y^2}{2} + \cdots + \frac{y^k}{k} \right)^{5} \frac{z^{k+1}}{(k+1)^4}\right].
\end{equation}

\section{Envoi: Two Research Items}

\textbf{Item 1:} \textit{Research project on Euler sums}.\index{Euler, L.}\index{Euler sum} Find a proof of the three experimentally obtained sums:
\begin{equation} \label{23.85a}
  s_h(3,2) = \frac{15}{2} \zeta(5) + \zeta(2) \zeta(3),
\end{equation}
\begin{equation} \label{23.86}
  s_h(4,3) = -\frac{109}{8} \zeta(7) + \frac{37}{2} \zeta(3) \zeta(4) -5 \zeta(2) \zeta(5),
\end{equation}
\begin{equation} \label{23.87}
  s_h(5,4) = \frac{890}{9}\zeta(9) +66\zeta(4)\zeta(5) -\frac{4295}{24}\zeta(2)\zeta(5) -5\zeta(3)^3 +\frac{265}{8} \zeta(2)\zeta(7).
\end{equation}

\textbf{Item 2:} \textit{Research project for limiting cases}. We know already that (\ref{23.85a}) to (\ref{23.87}) are cases of
\begin{equation} \label{23.88}
  s_h(m+1,m) = \sum_{k=1}^{\infty} \left( \frac{1}{1} + \frac{1}{2} + \cdots + \frac{1}{k} \right)^{m+1} \frac{1}{(k+1)^m}.
\end{equation}
Find the resulting infinite products as $y\rightarrow1$ and $z\rightarrow1$ for,
\begin{equation} \label{23.89}
\prod_{\substack{\gcd(j_1,j_2,j_3,k)=1 \\ j_1,j_2,j_3<k \\ j_1,j_2,j_3 \geq 1; k \geq 2}} \left(\frac{1}{1-y^{j_1+j_2+j_3}z^k} \right)^{\frac{1}{k}},
\end{equation}
\begin{equation} \label{23.90}
\prod_{\substack{\gcd(j_1,j_2,j_3,j_4,k)=1 \\ j_1,j_2,j_3,j_4<k \\ j_1,j_2,j_3,j_4\geq 1; k \geq 2}} \left(\frac{1}{1-y^{j_1+j_2+j_3+j_4}z^k} \right)^{\frac{1}{k}},
\end{equation}
\begin{equation} \label{23.91}
\prod_{\substack{\gcd(j_1,j_2,j_3,j_4,j_5,k)=1 \\ j_1,j_2,j_3,j_4,j_5<k \\ j_1,j_2,j_3,j_4,j_5\geq 1; k \geq 2}} \left(\frac{1}{1-y^{j_1+j_2+j_3+j_4+j_5}z^k} \right)^{\frac{1}{k}},
\end{equation}
and apply the three sums from the previous exercise to obtain three new closed form infinite products.

\end{document}